\documentclass{compositio}

\usepackage{amsmath, amssymb, amsfonts}
\usepackage{graphicx,bbm,color}
\newcommand{\RR}{\mathbb{R}}

\newcommand{\PP}{\mathbb{P}}
\newcommand{\EE}{\mathbb{E}}
\newcommand{\bfa}{\mbox{\boldmath$a$}}
\newcommand{\bfc}{\mbox{\boldmath$c$}}
\newcommand{\bfq}{\mbox{\boldmath$q$}}
\newcommand{\bfw}{\mbox{\boldmath$w$}}
\newcommand{\bbfw}{\mbox{\boldmath$\bar{w\hspace{.05pc}}\hspace{-.05pc}$}}
\newcommand{\tildebfc}{\mypp\tilde{\mynn\mbox{\boldmath$c$}}}
\newcommand{\bfdelta}{\mbox{\boldmath$\delta$}}
\newcommand{\bfOne}{\mathbbm{1}}
\newcommand{\bfone}{{\mathbf{1}}}

\newcommand{\fL}{{\mathfrak{h}}}

\newcommand{\calA}{\mathcal{A}}
\newcommand{\calF}{\mathcal{F}}
\newcommand{\calJ}{\mathcal{J}}
\newcommand{\calR}{\mathcal{R}}

\newcommand{\rme}{{\rm e}}
\newcommand{\rmd}{{\rm d}}

\newcommand{\myp}{\mbox{$\:\!$}}
\newcommand{\mypp}{\mbox{$\;\!$}}

\newcommand{\myn}{\mbox{$\;\!\!$}}
\newcommand{\mynn}{\mbox{$\:\!\!$}}

\newtheorem{theorem}{Theorem}[section]

\newtheorem{lemma}[theorem]{Lemma}

\theoremstyle{remark}
\newtheorem{remark}{Remark}[section]

\theoremstyle{definition}

\numberwithin{equation}{section}

\begin{document}

\title{Occupation time distributions for the telegraph process}
\author{Leonid Bogachev}
\email{L.V.Bogachev@leeds.ac.uk}
\address{Department of Statistics, University of Leeds, Leeds LS2 9JT, United Kingdom}

\author{Nikita Ratanov}
\email{nratanov@urosario.edu.co}
\address{Facultad de Economia, Universidad del Rosario, Cl.\,14, No.~4-69, Bogot\'a, Colombia}



\thanks{
The first author was partially supported by a Leverhulme Research
Fellowship and DFG Grant 436 RUS 113/722.
}

\begin{abstract}
For the one-dimensional telegraph process, we obtain explicit
distribution of the occupation time of the positive half-line. The
long-term limiting distribution is then derived when the initial
location of the process is in the range of sub-normal or normal
deviations from the origin; in the former case, the limit is given
by the arcsine law. These limit theorems are also extended to the
case of more general occupation-type functionals.
\end{abstract}

\maketitle

\section{Introduction}\label{sec:Intro}

Let $B=(B_t,\, t\ge 0)$ be a standard Brownian motion on $\RR$
starting from the origin ($B_0=0$), and consider the occupation time
functional
\begin{equation}\label{eq:h_T}
\fL_T:=\frac{1}{T}\int_0^T H(B_t)\,\rmd t,\qquad T>0,
\end{equation}
where $H(x)$ is the Heaviside unit step function (i.e., $H(x)=0$ for
$x\le0$ and $H(x)=1$ for $x>0$). That is to say, $\fL_T\in[0,1]$ is
the proportion of time spent by the Brownian motion $(B_t,\ 0\le
t\le T)$ on the positive half-line. It is well known  that the
probability distribution of the random variable $\fL_T$ does not
depend on $T$ (which is evident from the scaling property of the
Brownian motion and the fact that $H(\alpha x)\equiv H(x)$ for any
$\alpha>0$) and is given by the classic \textit{arcsine law},
\begin{equation}\label{Levy}
\PP\{\fL_T\le y\}=\frac{2}{\pi}\arcsin\sqrt{y},\qquad 0\leq y\leq 1,
\end{equation}
with the probability density
\begin{equation}\label{eq:arcsine}
p_{\mathrm{as}}(y):=\frac{1}{\pi\sqrt{y\myp(1-y)}}\,,\qquad 0<y<1.
\end{equation}

The beautiful formula (\ref{Levy}) dates back about 70 years to
P.~L\'evy \cite[Th\'eor\`eme~3, pp.\:301--302]{Levy}, who has also
proved that the arcsine law (\ref{Levy}) is the limit distribution
for the relative frequency of positive sums among consecutive
partial sums of independent symmetric Bernoulli ($0\text{--}1$)
random variables \cite[Corollaire 2, p.\;303]{Levy}. Using the
invariance principle, the latter result was extended by P.~Erd\H{o}s
and M.~Kac \cite{EK} to the case of sums of arbitrary i.i.d.\ random
variables with zero
mean and unit variance (cf.\ \cite[
Theorem 4.3.19, p.\;236]{Stroock}). More recently, R.~Khasminskii
\cite{Khas} obtained the limit distribution, as $T\to\infty$, of
more general functionals of the form
$$
\fL_T (x;f):=\frac{1}{T}\int_0^T f(x+X_t)\,\rmd t,
$$
where $(X_t,\,t\ge0)$ is a diffusion process on $\mathbb{R}$
($X_0=0$) with generator $L=-a(x)\,\rmd^2/\rmd x^2$, and
$f:\RR\to\RR$ is a probing function from a suitable class. In
particular, the results of \cite{Khas} imply that if
$\lim_{x\to\pm\infty}a(x)=a_0>0$ and $f$ is a bounded piecewise
continuous function such that
\begin{equation}\label{Cesaro}
\lim_{x\to\pm\infty} \frac{1}{x}\int_0^x f(u)\,\rmd{u} =
f_\pm\,,\qquad f_+\ne f_-\,,
\end{equation}
then the distribution of the random variable
$(\fL_T(x;f)-f_{-})/(f_{+}-f_{-})$ converges weakly, as
$T\to\infty$, to the arcsine law (\ref{Levy}).

In the present paper, we obtain similar results for the so-called
\emph{telegraph process} defined by
\begin{equation}\label{VX}
X_t:= V_0 \int_0^t(-1)^{N_u}\,\rmd u,\qquad t\ge 0,
\end{equation}
where $(N_t,\,t\ge 0)$ is a homogeneous Poisson process (with rate
$\lambda>0$), $V_0$ is a random variable with equiprobable values
$\pm c$ independent of the process $N_t$, and $c>0$ is a parameter
(see \cite{Gold,Kac1,Pinsky}). That is, $X_t$ is the position at
time $t\ge0$ of a particle starting at $t=0$ from the origin and
moving on the line with alternating velocities $\pm c$, reversing
the direction of motion at each jump instant of the Poisson process
$N_t$; the initial (random) direction is decided by the sign of
$V_0$. Note that the process $X_t$ itself is non-Markovian, however
if $V_t= \rmd X_t/\rmd{t}=(-1)^{N_t}\myp V_0$ is the corresponding
velocity process, then the joint process $(X_t,V_t)$ is Markov on
the state space $\RR\times\{-c,+c\}$ (see \cite[\S\myp12.1,
p.\;469]{EtK}). We shall also consider the conditional telegraph
processes obtained from $X_t$ by conditioning on $V_0$,
\begin{equation}\label{VX-pm}
X_t^\pm:=\pm\mypp c\int_0^t(-1)^{N_u}\,\rmd u,\qquad t\ge 0,
\end{equation}
where the choice of the $+$ or $-$ sign determines the initial
direction of motion.

\begin{remark}
Here and throughout the paper, we adopt a notational convention that
any formula involving the $\pm$ and $\mp$ signs combines the two
cases corresponding to the choice of either the upper or lower sign,
respectively.
\end{remark}

\begin{remark}
The telegraph process is the simplest example of so-called
\textit{random evolutions} (see, e.g., \cite[Ch.\,12]{EtK} and
\cite[Ch.\;2]{Pinsky}).
\end{remark}

The model of non-interacting particles moving in one dimension with
alternating velocities (updated at random on a discrete time grid)
was first introduced in 1922 by G.~I.~Taylor \cite{Taylor} in an
attempt to describe turbulent diffusion; later on (around
1938--1939) it was studied at length by S.~Goldstein \cite{Gold} in
connection with a certain hyperbolic partial differential equation
(called the \textit{telegraph}, or \textit{damped wave}
\textit{equation}, see (\ref{eq:te}) below) describing the
spatio-temporal dynamics of the potential in a transmitting cable
(without leakage) \cite{Webster}. In his 1956 lecture notes, M.~Kac
(see \cite{Kac1}) considered a continuous-time version of the
telegraph model. Since then, the telegraph process and its many
generalizations have been studied in great detail (see, e.g.,
\cite{Or90,Pinsky,Or95,MPRF,Weiss}), with numerous applications in
physics \cite{Weiss},
biology \cite{Hadeler-rev,HH}, ecology \cite{OL} and, more recently,
financial market modelling \cite{QF,RM08} (see also further
bibliography in these papers).

An efficient conventional approach to the analytical study of the
telegraph process, analogous to that for diffusion processes, is
based on pursuing a fundamental link relating various expected
values of the process with initial value and/or boundary value
problems for certain partial differential equations (see, e.g.,
\cite{Gold,Or90,Or95,R3,MPRF,Ra06}). In particular, Kac \cite{Kac1}
has shown that, for any bounded continuously differentiable function
$g_0:\RR\to\RR$, the functions
\begin{equation*}
v^\pm(x,t):=\EE\bigl[g_0(x+X^\pm_t)\bigr],\qquad x\in\RR,\ t\ge0,
\end{equation*}
satisfy the set of partial differential equations
\begin{equation*}
\frac{\partial v^\pm(x,t)}{\partial t}\mp c\,\frac{\partial
v^\pm(x,t)}{\partial
x}=\mp\lambda\left(v^{+}(x,t)-v^{-}(x,t)\right), \qquad t>0,
\end{equation*}
with the initial conditions
\begin{equation*}
v^\pm(x,0)=g_0(x), \qquad x\in\RR.
\end{equation*}
These equations can be easily combined (see details in \cite{Kac1}
or \cite[\S\myp12.1, p.\;470]{EtK}) to show that the function
\begin{equation*}
v(x,t):=\EE\left[g_0(x+X_t)\right]=\frac12\,v^-(x,t)+\frac12\,v^+(x,t)
\end{equation*}
satisfies the telegraph (or telegrapher's) equation (see, e.g.,
\cite[\S\myp15]{Webster})
\begin{equation}\label{eq:te}
\frac{\partial^2 v}{\partial t^2}+2\lambda \frac{\partial
v}{\partial t}=c^2 \frac{\partial^2 v}{\partial x^2}\,
\end{equation}
with the initial conditions
\begin{equation}\label{ICH}
v(x, 0)=g_0(x),\qquad \frac{\partial v}{\partial t}(x, 0)=0\myp.
\end{equation}

\begin{remark}
The telegraph equation (\ref{eq:te}) first appeared more than 150
years ago in work by W.~Thomson (Lord Kelvin) on the transatlantic
cable \cite{Kelvin}.
\end{remark}

The (unique) solution of the Cauchy problem
(\ref{eq:te})\myp--\myp(\ref{ICH}) can be written explicitly (see,
e.g., \cite[\S\S\,46, 74]{Webster} or \cite[\S\,0.4]{Pinsky})
as
\begin{equation}\label{sol:te}
\begin{aligned}
v(x, t)&=\frac{1}{2}\,\rme^{-\lambda t}\bigl(
g_0(x+ct)+g_0(x-ct)\bigr)\\
&\quad +\frac{1}{2}\,\rme^{-\lambda t}\int_{-t}^tg_0(x+cu)\left(
\lambda I_0\bigl(\lambda\sqrt{t^2-u^2}\,\bigr)+ \frac{\lambda\myp
t}{\sqrt{t^2-u^2}}\,I_1\bigl(\lambda\sqrt{t^2-u^2}\bigr)
\right)\rmd{u}\myp,
\end{aligned}
\end{equation}
where
\begin{equation*}
I_0(z):=\sum_{n=0}^\infty \frac{(z/2)^{2n}}{(n!)^{2}}\ \ \quad
\text{and} \ \ \quad I_1(z):=I_0'(z)=\frac{z}{2}\sum_{n=0}^\infty
\frac{(z/2)^{2n}}{n!\,(n+1)!}\qquad (z\in\RR)
\end{equation*}
are the modified Bessel functions of the first kind (of orders $0$
and $1$), respectively) \cite[
9.6.12, p.\;375; 9.6.27,
p.\;376]{AS}.


It is well known that, under a suitable scaling, the telegraph
process satisfies a functional central limit theorem.
\begin{theorem}\label{Kac}
Assume that\/ $\lambda,\myp{}c\to +\infty$ in such a way that
$c^2\myn/\lambda\to 1$. Then the distribution of the telegraph
processes $(X^\pm_t,\,t\ge 0)$ converges weakly in $C[0,\infty)$ to
the distribution of a standard Brownian motion $(B_t,\,t\ge0)$. The
same is true for the unconditional telegraph process
$(X_t,\,t\ge0)$.
\end{theorem}

As was observed by Kac \cite[p.\;501]{Kac1}, this result formally
follows from the telegraph equation (\ref{eq:te}),
which in the limit $\lambda,\mypp{}c\to +\infty$,
\,$c^2\myn/\lambda\to 1$ yields the diffusion (heat) equation
\begin{equation*}
\frac{\partial v}{\partial t}=\frac{1}{2}\, \frac{\partial^2
v}{\partial x^2}\,,
\end{equation*}
associated with the standard Brownian motion $B_t$. A rigorous proof
of Theorem \ref{Kac}, along with some extensions, can be found in
\cite[\S\myp12.1, p.\,471]{EtK} and \cite[Theorem 5.1]{MPRF}.

Our main goal in the present paper is to analyze the distribution of
the occupation time of the telegraph process $X_t$ and, in
particular, to obtain a limit distribution, as $T\to\infty$, of the
occupation-type functionals of the form
$\eta_T(x;f):=T^{-1}\!\int_0^T f(x+X_t)\,\rmd t$ for a suitable
class of probing functions~$f$. In particular, we prove that the
limit distribution is given by L\'evy's arcsine law providing that
the starting point $x$ is in the range of subnormal deviation from
the origin (i.e., $x=o(\sqrt{T}\,)$). For technical simplicity, we
impose a stronger condition on the asymptotics of $f$ at
$\pm\infty$, assuming that the corresponding limits $f_{\pm}$ exist.

The rest of the paper is organized as follows. In Section
\ref{sec:Results} we state the main results of this work (Theorems
\ref{th1}\myp--\myp\ref{th4}), which are then proved in Sections
\ref{sec:th1}\myp--\myp\ref{sec:th4}, respectively. Section
\ref{sec:FK} contains a suitable version of the Feynman--Kac
formula, with applications to the Laplace transforms for the
occupation-type functionals under study, which is instrumental for
our techniques. We finish in Section \ref{sec:Concl} with concluding
remarks and some conjectures, which are illustrated by the results
of computer simulations. Appendices \ref{sec:A1} and \ref{sec:A2}
contain alternative (probabilistic) proofs of Theorems \ref{th2} and
\ref{th3}, respectively.

\section{Statement of the main results}\label{sec:Results}

For $T>0$, \,$x\in\RR$, consider the following occupation time
random variables
\begin{equation}\label{OT0}
\eta_T(x):=\frac{1}{T}\int_0^{T} H(x+X_t)\,\rmd t,\qquad
\eta^\pm_T(x):=\frac{1}{T}\int_0^{T} H(x+X^\pm_t)\,\rmd t,
\end{equation}
where $H(x)=\bfOne_{(0,\infty)}(x)$ is the Heaviside step function
and $X_t$, $X_t^\pm$ are the telegraph processes introduced above
(see (\ref{VX}), (\ref{VX-pm})). Note that the total time spent by
the processes $(x+X^\pm_t,\,0\le t\le T)$ at the origin almost
surely (a.s.) equals zero, since by Fubini's theorem we have
\begin{equation}\label{eq:meas=0}
\EE\int_0^{T}\!\bfOne_{\{0\}}(x+X^\pm_t)\,\rmd t=\int_0^{T}
\PP\{X^\pm_t=-x\} \,\rmd t=0\mypp.
\end{equation}
Hence, the complementary quantity $1-\eta^\pm_T(x)$ a.s.\ equals the
proportion of time spent by the processes $(x+X^\pm_t,\,0\le t\le
T)$ on the negative side of the axis,
\begin{equation*}
1-\eta^\pm_T(x)=\frac{1}{T}\int_0^{T}
\bfOne_{(-\infty,0)}(x+X^\pm_t)\,\rmd t\qquad\text{(a.s.)},
\end{equation*}
and by symmetry (with respect to simultaneous transformations
$x\mapsto -x$, \,$\pm\mapsto\mp$) it follows that
\begin{equation}\label{eq:=d}
\eta^\pm_T(x)\stackrel{d}{=}1-\eta^\mp_T(-x),\qquad x\in\RR\mypp.
\end{equation}

Let us consider the function $\varphi_T(t)$ ($t\ge0$) defined by
\begin{equation}\label{eq:phi_T}
\varphi_T(t):=\frac{1}{4\pi\lambda T}\int_0^t\frac{1-\rme^{-2\lambda
Tu}}{u^{3/2}\sqrt{t-u}}\,\rmd u \ \quad (t> 0),\qquad
\varphi_T(0):=\frac12\,.
\end{equation}
After the substitution $u=ty$, we have in the limit as
$t\downarrow0$,
\begin{equation}\label{eq:phi_T1}
\varphi_T(t)=\frac{1}{4\pi\lambda Tt}\int_0^1\frac{1-\rme^{-2\lambda
Tty}}{y^{3/2}\sqrt{1-y}}\,\rmd{y}\to
\frac{1}{2\pi}\int_0^1\frac{1}{\sqrt{y(1-y)}}\,\rmd{y}=\frac{1}{2}
\end{equation}
(see (\ref{eq:arcsine})), and so $\varphi_T(\cdot)$ is continuous at
zero (and hence everywhere on $[0,\infty)$). Note the following
useful scaling relation, which easily follows from the
representation of $\varphi$ given by (\ref{eq:phi_T1}):
\begin{equation}\label{eq:phi-scaling}
\varphi_{\alpha T}(t)=\varphi_{T}(\alpha t),\qquad t\ge 0, \quad
\alpha>0\myp.
\end{equation}
Let us also set
\begin{equation}\label{eq:psi_T}
\psi_T(y):=2\lambda T\varphi_T(y) \mypp\varphi_T(1-y),\qquad 0\le
y\le 1\myp.
\end{equation}

We are now ready to state our first result.
\begin{theorem}\label{th1}
The random variables\/ $\eta^\pm_T(0)$ defined in
\textup{(\ref{OT0})} have the distribution
\begin{equation}
\label{eta0a}
\PP\bigl\{\eta^{\pm}_T(0)\in{\rmd}y\bigr\}=2\myp\varphi_T(1)\,\delta_{x^{\pm}}
(\rmd{y})+ \psi_T(y)\,\rmd{y}\myp,\qquad 0\le y\le 1,
\end{equation}
where $\delta_x$ is the Dirac measure \textup{(}of unit
mass\textup{)} at point $x$\textup{,} with
$x^{-}\mynn=0$\myp\textup{,} $x^{+}\mynn=1$\myp. Furthermore, the
distribution of $\eta_T(0)$ \textup{(}see
\textup{(\ref{OT0})}\textup{)} is given by the formula
\begin{equation}\label{eta00-}
\PP\bigl\{\eta_T(0)\in
\rmd{y}\bigr\}=\varphi_T(1)\mypp\delta_{0}(\rmd{y})+
\varphi_T(1)\mypp\delta_{1}(\rmd{y})+\psi_T(y)\,\rmd{y}\myp,\qquad
0\leq y\le 1.
\end{equation}
\end{theorem}

In other words, the distribution of $\eta_T^-(0)$, $\eta_T^+(0)$ has
a discrete part with atom of mass $2\myp\varphi_T(1)$ at point $0$
or $1$, respectively, and an absolutely continuous part with the
density $\psi_T$ defined by (\ref{eq:psi_T}). Similarly, the
distribution of $\eta_T(0)$ has atoms at points $0$ and $1$, both of
mass $\varphi_T(1)$, and an absolutely continuous part with the
density $\psi_T$ as above.

\begin{remark}
The $\pm$-duality in (\ref{eta0a}) becomes clear from relation
(\ref{eq:=d}) (with $x=0$) and the symmetry property
$\psi_T(y)\equiv\psi_T(1-y)$ (see (\ref{eq:psi_T})).
\end{remark}

\begin{remark}
Using an integral formula (see
\cite[9.6.16, p.\;376]{AS}) for the modified Bessel function $I_0$,
it is easy to check that the function $\varphi_T$ defined by
(\ref{eq:phi_T}) admits another representation,
\[
\varphi_T(t)=\frac{1}{2\lambda Tt}\int_0^{\lambda Tt}\rme^{-y}\mypp
I_0(y)\,\rmd{y}\myp,\qquad t>0\myp,
\]
which is further evaluated (see \cite[11.3.12, p.\;483]{AS}) to
yield $\varphi_T(t)=\frac{1}{2}\mypp\rme^{-\lambda
Tt}\bigl(I_0(\lambda Tt)+I_1(\lambda Tt)\bigr)$. Thus, the
distribution of $\eta_T^\pm(0)$ and $\eta_T(0)$ can be expressed
through the modified Bessel functions $I_0$ and $I_1$, as well as
the distribution of the telegraph process (cf.\ (\ref{sol:te})).
\end{remark}

In the next theorem, we give explicit integral formulas for the
distribution of $\eta_T^{\pm}(x)$ in the case $x\ne0$. For
simplicity, we only present the answer for $x<0$, the case $x>0$
readily following in view of the duality relation (\ref{eq:=d}).

\begin{theorem}\label{th2}
Assume that $x<0$ and set $T_0:=|x|/c$. Then\textup{,} for any
$T>0$\myp\textup{,} the random variables $\eta^{\pm}_T(x)$ defined
in \textup{(\ref{OT0})} have the following distribution\textup{:}

\textup{(a)} \,if\/ $T\le T_0$ then
$\PP\{\eta^{\pm}_T(x)=0\}=1$\myp\textup{;}

\textup{(b)} \,if\/ $T>T_0$ then\textup{,} for\/ $0\le y\le
1-T_0/T$\textup{,}
\begin{equation}\label{eq:Y-}
\PP\{\myp{}\eta_T^{\pm}(x)\in \rmd{y}\}= \left(\int_T^\infty
Q^{\pm}_{-x}(u)\,\rmd{u}\right)\delta_{0}(\rmd{y})+\mu_T^{\pm}(\rmd{y})
+\varPsi_{x}^{\pm}(y,T)\mypp\rmd{y}\myp,
\end{equation}
where $\mu_T^{-}(\rmd{y}):= 0$ and
\begin{gather}\label{eq:mu}
\mu_T^{+}(\rmd{y}):=2\mypp\rme^{-\lambda T_0}
\varphi_{T}(1-T_0/T)\,\delta_{1-T_0/T}(\rmd{y})+\rme^{-\lambda T_0}
\,\psi_{T-T_0}\left(\frac{y}{1-T_0/T}\right)\frac{\rmd{y}}{1-T_0/T}\,,\\
\label{eq:varPsi} \varPsi_{x}^{\pm}(y,T):=2\myp T\myp
Q^{\pm}_{-x}((1-y)T)\,\varphi_{T}(y)+\int_{T_0}^{(1-y)\myp T}
\!Q^{\pm}_{-x}(u)\,\psi_{T-u}
\left(\frac{y}{1-u/T}\right)\frac{\rmd{u}}{1-u/T}\,,
\end{gather}
with $\varphi_T$ and $\psi_T$ given by \textup{(\ref{eq:phi_T})} and
\textup{(\ref{eq:psi_T})}\textup{,} respectively\textup{,} and the
functions $Q_{-x}^{\pm}(u)$ $(-x>0)$ defined for all
$u\in[T_0,\infty)$ by
\begin{align}\label{eq:Q+}
Q^+_{-x}(u)&:=\frac{\lambda T_0\,\rme^{-\lambda
u}}{\sqrt{u^2-T_0^2}}\,I_1\bigl(\lambda
{\textstyle\sqrt{u^2-T_0^2}}\,\bigr)\myp,\\
\label{eq:Q-} Q^-_{-x}(u)&:=\lambda\mypp\rme^{-\lambda
u}\,I_0\bigl(\lambda
{\textstyle\sqrt{u^2-T_0^2}}\,\bigr)-\frac{\lambda
(u-T_0)\,\rme^{-\lambda u}}{\sqrt{u^2-T_0^2}}\,I_1\bigl(\lambda
{\textstyle\sqrt{u^2-T_0^2}}\,\bigr)\myp.
\end{align}
\end{theorem}

For the next theorem, we need a few notations. For $a>0$, consider
the function
\begin{equation}\label{eq:q}
q_{a}(t):=\frac{a}{\sqrt{2\pi
t^3}}\,\exp\left(-\frac{a^2}{2t}\right), \qquad t>0,
\end{equation}
with Laplace transform (see \cite[29.3.82, p.\,1026]{AS})
\begin{equation}\label{eq:4.9}
\int_0^\infty \rme^{-s\myp t}\mypp
q_{a}(t)\,\rmd{t}=\rme^{-a\sqrt{2s}}\,,\qquad s\ge0.
\end{equation}
Let $Y_a$ ($a\ge0$) be a family of random variables with values in
$[0,1]$, such that $Y_0$ has the arcsine distribution (\ref{Levy}),
with the density $p_{\mathrm{as}}$ (see (\ref{eq:arcsine})), while
for $a>0$ the distribution of $Y_a$ is given by
\begin{equation}\label{eq:Y0}
\PP\{Y_a\in \rmd{y}\}= m_a\mypp\delta_{0}(\rmd{y})+ f_a(y)\,\rmd{y},
\end{equation}
where
\begin{align}
\label{eq:m_a} m_a:=\int_{1}^\infty
q_{a}(u)\,\rmd{u}&=\frac{2}{\sqrt{2\pi}}\int_0^{a}\rme^{-y^2/2}\,\rmd{y}\myp,\\
\label{eq:f_a} f_a(y):=\int_0^{1-y}
\frac{q_{a}(u)}{1-u}\,p_{\mathrm{as}}\!\left(\frac{y}{1-u}\right)\rmd{u}
&=\frac{a}{\sqrt{2\pi^3y}}\int_0^{1-y}
\!\frac{\rme^{-a^2\mynn/(2u)}}{u^{3/2}\,\sqrt{1-y-u}}\:\rmd{u}\mypp.
\end{align}

\begin{remark}
It is easy to verify, either from (\ref{eq:q}) or using the Laplace
transform (\ref{eq:4.9}), that $q_a\stackrel{w^*}{\to} \delta_0$ as
$a\to0+$, where $\delta_0(\cdot)$ is the Dirac delta function and
\,$\stackrel{w^*}{\to}$\, denotes weak-$^*$ convergence of
generalized functions; hence $m_a\to0$ (which can also be seen
directly from the right-hand side of (\ref{eq:m_a})) and
$f_a\stackrel{w^*}{\to} p_{\mathrm{as}}$ (see the first part of
formula (\ref{eq:f_a})). That is, $Y_a\stackrel{d\,}{\to}Y_0$ as
$a\to0+$, and so the distribution of $Y_a$ is continuous in
parameter $a\in[0,\infty)$.
\end{remark}

\begin{theorem}\label{th3}
Suppose that the initial position $X^\pm_0=x$\textup{,} as well as
the parameters $c$ and $\lambda$\textup{,} may depend on $T$ in such
a way that $\lambda T\to\infty$ and $(c^2T\myn/\lambda)^{-1/2}\,x\to
a\in\RR$ as $T\to\infty$. Then\textup{,} as $T\to\infty$\textup{,}
\begin{equation}\label{eq:T->infty}
\eta_T(x),\, \eta_T^\pm(x)\stackrel{d}{\longrightarrow}
\left\{\begin{array}{ll} Y_{-a}\,,& a\le 0\myp,\\[.2pc]
1-Y_a\mypp,\quad& a\ge0\myp.
\end{array}
\right.
\end{equation}
In particular\textup{,} for $a=0$ the limit is given by the arcsine
distribution \textup{(\ref{Levy})}.
\end{theorem}

To order to generalize these results in the spirit of \cite{Khas},
let $f:\RR\to\RR$ be a bounded, piecewise continuous function (i.e.,
continuous on $\RR$ outside a finite set $D_{\myn f}$, where it has
finite left and right limits), such that, for some finite constants
$f_+\neq f_-$\,,
\begin{equation}\label{assump}
\lim_{x\to-\infty} f(x)=f_{-}\,,\qquad \lim_{x\to+\infty}
f(x)=f_{+}\,.
\end{equation}
Consider the random variables
\begin{equation}\label{etaf}
\eta^\pm_T(x;f):=\frac{1}{(f_{+}-f_{-})\, T}
\int_0^T\bigl(f(x+X^\pm_t)-f_-\bigr)\,\rmd{t},\qquad x\in\RR.
\end{equation}
Clearly, by a linear transformation of the function,
\begin{equation}\label{eq:std}
f(x)\mapsto\tilde f(x):=\frac{f(x)-f_-}{f_+-f_-}\,,\qquad x\in\RR,
\end{equation}
we may and will
assume without loss of generality that $f_-=0$, \,$f_+=1$, so that
(\ref{etaf}) is reduced to
\begin{equation}\label{etaf'}
\eta^\pm_T(x;f):=\frac{1}{T} \int_0^T \!f(x+X^\pm_t)\,\rmd t.
\end{equation}

\begin{theorem}\label{th4}
Let the function $f$ satisfy the above conditions including
assumption \textup{(\ref{assump})} with $f_{-}=0$, $f_{+}=1$.
Suppose that the hypotheses of Theorem \textup{\ref{th3}} are
satisfied\myp\textup{,} and assume in addition that\/
$c^2T\myn/\lambda\to\infty$ as $T\to\infty$. Then the distribution
of\/ $\eta^\pm_T(x;f)$ converges weakly\myp\textup{,} as
$T\to\infty$\textup{,} to the law determined by the right-hand side
of \textup{(\ref{eq:T->infty})}.
\end{theorem}

\begin{remark}
Theorem \ref{th3} may be inferred from the diffusion approximation
(Theorem \ref{Kac}) of the telegraph process (see an alternative
proof in the Appendix \ref{sec:A2}). However, we will give a direct
proof of Theorem \ref{th3}, which may be of interest in its own
right and will also be instrumental for laying out the necessary
techniques for the proof of Theorem \ref{th4}, where the ``diffusion
approximation'' trick does not seem to be readily applicable.
\end{remark}

\section{The Feynman--Kac formula and applications}\label{sec:FK}

Let us recall the Feynman--Kac formula for the telegraph processes.

\begin{theorem}\label{th:FK}
Let $(X^\pm_t,\,t\ge0)$ be the telegraph processes
\textup{(\ref{VX-pm})}. Suppose that $g_0$ and $g$ are bounded
functions on $\RR$ such that $g_0\in C^1(\RR)$ and $g$ is piecewise
continuous\textup{,} i.e.\textup{,} $g\in C(\RR\setminus\myn
D_g)$\textup{,} where $D_g$ is a finite set\textup{,} and
moreover\textup{,} $f$ has finite left and right limits at the
points of $D_g$. Then the functions
\begin{equation}\label{eq:v-pm}
v^\pm(x, t):=\EE\left[g_0(x+X^\pm_t)\exp\left\{\int_0^t
g(x+X^\pm_u)\,\rmd u\right\}\right],\qquad x\in\RR,\ \ t\ge0,
\end{equation}
for all $(x,t)\in\RR\times\RR_+$ such that $x\pm ct\notin D_g$
satisfy the set of partial differential equations
\begin{equation*}
\frac{\partial v^{\pm}(x, t)}{\partial t}\mp c\,\frac{\partial
v^{\pm}(x, t)}{\partial x}=\mp\lambda\left(v^{+}(x, t)-
v^{-}(x,t)\right)+g(x)\mypp v^{\pm}(x, t),
\end{equation*}
with the initial conditions
\begin{equation*}
v^\pm(x,0)=g_0(x),\qquad x\in \RR.
\end{equation*}
\end{theorem}

This theorem is proved (see details in \cite{Ra06}) similarly to the
analogous result for diffusion processes (cf., e.g.,
\cite[\S\mypp2.6]{Ito}). An alternative probabilistic representation
for the solution of a deterministic telegraph-like equation is
developed in \cite{DMT}.

Let $\eta_T^\pm(x)$ be defined by (\ref{OT0}). For $\beta\in\RR$,
set
\begin{equation}\label{eq:v}
v^\pm_T(\xi,t):= \EE\bigl[\rme^{-\beta
t\mypp\eta^\pm_{Tt}(cT\xi)}\bigr],\qquad \xi\in\RR,\ \ t\ge0,
\end{equation}
or more explicitly (cf.\ (\ref{eq:v-pm}))
\begin{equation}\label{eq:v1}
v^\pm_T(\xi,t)= \EE\left[\exp\left\{\frac{-\beta}{T}\int_0^{Tt}
H(cT\xi+X^\pm_u)\,\rmd u\right\}\right],\qquad \xi\in\RR,\ \ t\ge0.
\end{equation}
Since $H(\cdot)$ is a bounded function, the expectation in
(\ref{eq:v1}) is finite for all $\beta\in\RR$.

Let us record some simple properties of the function $v^\pm_T$\myp.
\begin{lemma}\label{lm:pm}
For each $\beta\in\RR$ and any $T>0$\myp\textup{,} the functions
$v^\pm_T(\xi,t)$ are continuous on $\RR\times\RR_+$ and
\begin{equation}\label{eq:pm}
\lim_{\xi\to-\infty} v^\pm_T(\xi,t)=1\myp,\qquad
\lim_{\xi\to+\infty}v^\pm_T(\xi,t)=\rme^{-\beta t}\myp.
\end{equation}
\end{lemma}

\proof Continuity in $t\in\RR_+$ is obvious. As mentioned above (see
(\ref{eq:meas=0})), for any $\xi_0\in\RR$ we have a.s.\ that
$cT\xi_0+X^\pm_u\ne0$ for all $u\in[0,Tt]$ except on a (random) set
of Lebesgue measure zero. Since the function $H$ is continuous
outside zero, this implies that, for such $u$,
$H(cT\xi+X^\pm_u)\stackrel{\text{\!a.s.\,}}{\longrightarrow}
H(cT\xi_0+X^\pm_u)$ as $\xi\to\xi_0$ and hence, by Lebesgue's
dominated convergence theorem, $\int_0^{Tt}
H(cT\xi+X^\pm_u)\,\rmd{u}\stackrel{\text{\!a.s.\,}}{\longrightarrow}
\int_0^{Tt} H(cT\xi_0+X^\pm_u)\,\rmd{u}$ as $\xi\to\xi_0$. The
continuity of $v_T^\pm(\cdot,t)$ at point $\xi_0$ now follows by
Lebesgue's dominated convergence theorem applied to the expectation
(\ref{eq:v1}), since everything is bounded (for a fixed $t$).

To prove (\ref{eq:pm}), note that, for $T>0$ and each $u\ge0$, we
have $cT \xi+X_u^{\pm}\stackrel{\text{\!a.s.\,}}{\longrightarrow}
\pm\infty$ as $\xi\to\pm\infty$. Since $H$ is bounded on
$\mathbb{R}$, the claim now follows by dominated convergence.
\endproof

From the definition (\ref{eq:v1}), it is clear that if $\beta\ge0$
then, for each $\xi\in\RR$, the functions $v_T^\pm(\xi,\cdot)$ are
bounded on $[0,\infty)$, so we can define the Laplace transform
\begin{equation}\label{eq:w}
w_T^\pm(\xi,s):=\int_0^\infty\rme^{-st}\mypp
v_T^\pm(\xi,t)\,\rmd{t}\qquad (s>0).
\end{equation}

\begin{lemma}\label{lm:tilde}
Set\/ $\tilde{s}:=s+\beta$. For any fixed $s>0$\myp\textup{,} the
functions $w_T^\pm=w_T^\pm(\xi,s)$ defined by \textup{(\ref{eq:w})}
are continuous in\/ $\xi\in\RR$ and satisfy the following set of
differential equations
\begin{equation}\label{sysw_comp}
\frac{\partial\myp w_T^{\pm}}{\partial \xi} =\lambda
T\bigl(w_T^{+}-w_T^{-}\bigr)\pm\bigl(s+H(cT\xi)\bigr)\mypp
w_T^{\pm}\mp 1,\qquad \xi\ne0.
\end{equation}
Moreover\myp\textup{,}
\begin{equation}\label{wc}
\lim_{\xi\to-\infty} w_T^\pm(-\xi,s)=s^{-1},\qquad
\lim_{\xi\to+\infty}w_T^\pm(\xi,s)=\tilde{s}^{\myp-1}.
\end{equation}
\end{lemma}

\proof The continuity of the functions $w_T^\pm(\xi,s)$ in $\xi$
follows from the definition (\ref{eq:w}) and the first part of Lemma
\ref{lm:pm}. Further, applying Theorem \ref{th:FK} (with
$g_0(x)\equiv 1$ and $g(x)=-\beta\myp T^{-1} H(x)$), we see that the
functions $v_T^\pm=v_T^\pm(\xi,t)$ defined by (\ref{eq:v}) satisfy
the initial value problem
\begin{gather}
\label{FKs1} \mp\frac{\partial v_T^{\pm}}{\partial t}+
\frac{\partial v_T^{\pm}}{\partial\xi}=\lambda T\left(v_T^{+}-
v_T^{-}\right)\pm\beta H(cT\xi)\mypp v_T^{\pm},\qquad t>0,\quad
\xi\pm t\ne 0,\\
\label{FKs-IC}
v_T^\pm(\xi,0)=1,\qquad \xi\in\RR.
\end{gather}
Integrating by parts and using the initial condition (\ref{FKs-IC}),
we have
\begin{equation}\label{eq:-1}
\int_0^\infty \rme^{-st}\,\frac{\partial\myp
v_T^\pm(\xi,t)}{\partial t}\,\rmd t= -v_T^\pm(\xi,0)+s\int_0^\infty
\rme^{-st}\,v_T^\pm(\xi,t)\,\rmd t =-1+sw_T^\pm(\xi,s).
\end{equation}
Applying the Laplace transform (with respect to $t$) to equation
(\ref{FKs1}) and taking into account (\ref{eq:-1}), we immediately
obtain the differential equation (\ref{sysw_comp}). Finally, the
boundary conditions (\ref{wc}) readily follow from (\ref{eq:pm}) by
Lebesgue's dominated convergence theorem applied to (\ref{eq:w}).
\endproof

Let us also make similar preparations for the random variables
$\eta_T^\pm(x;f)$ defined in (\ref{etaf}). As explained in Section
\ref{sec:Results} (see (\ref{eq:std})), without loss of generality
this definition can be simplified to the form (\ref{etaf'}).
Consider the functions (cf.\ (\ref{eq:v}))
\begin{equation*}
v_T^\pm(\xi,t;f):=\EE\bigl[\exp\left(-\beta
t\eta^\pm_{Tt}(cT\xi;f)\right)\bigr],\qquad \xi\in\RR,\ \ t\ge0,
\end{equation*}
and the corresponding Laplace transform
\begin{equation*}
w_T^\pm(\xi,s;f):=\int_0^\infty\rme^{-st}\mypp
v_T^\pm(\xi,t;f)\,\rmd t,\qquad s>0.
\end{equation*}
Then, again applying Theorem \ref{th:FK} (with $g_0(x)\equiv 1$ and
$g(x)=-\beta\myp T^{-1}f(x)$), similarly to Lemmas \ref{lm:pm} and
\ref{lm:tilde} one can show that $w_T^{\pm}= w_T^\pm(\xi,s;f)$, for
each $s>0$, is a continuous bounded function of\/ $\xi\in\RR$,
satisfying the differential equation (cf.\ (\ref{sysw_comp}))
\begin{equation}\label{sysw}
\frac{\partial\myp w_T^\pm}{\partial \xi} =\lambda
T\bigl(w_T^{+}-w_T^{-}\bigr)\pm\bigl(s+\beta f(cT\xi)\bigr)\mypp
w_T^{\pm}\mp 1,\qquad \xi\in \RR\setminus D_{\myn f}\myp,
\end{equation}
with the same boundary conditions at $\pm\infty$ as (\ref{wc}),
\begin{equation}\label{wc*}
\lim_{\xi\to-\infty}w_T^\pm(\xi,s;f)=s^{-1},\qquad
\lim_{\xi\to\infty}w_T^\pm(\xi,s;f)=\tilde{s}^{\myp-1}.
\end{equation}

\section{Proof of Theorem \ref{th1}}\label{sec:th1}
In what follows, the prime $'$ denotes the transposition of vectors.
Introducing the vector notations
\begin{equation*}
\bfw_T(\xi,s):=(w_T^+(\xi,s), w_T^-(\xi,s))', \qquad \bfone:=(1,1)',
\qquad \tilde{\bfone}:=(1,-1)',
\end{equation*}
we can write down equations (\ref{sysw_comp}) and (\ref{wc}) in the
matrix form,
\begin{gather}\label{sysw-m0}
\frac{\partial\myp\bfw_T(\xi,s)}{\partial \xi}=\calA_T(\xi,s)\mypp
\bfw_T(\xi,s)-\tilde{\bfone}\qquad (\xi\neq 0),\\
\label{limw-1} \lim_{\xi\to-\infty}\bfw_T(\xi,s)
=s^{-1}\myp\bfone\myp,\qquad \lim_{\xi\to+\infty}\bfw_T(\xi,s)
=\tilde{s}^{\myp-1}\myp\bfone\myp,
\end{gather}
where $\tilde{s}=s+\beta$ (see Lemma \ref{lm:tilde}) and
\begin{equation}\label{eq:A0}
\calA_T(\xi,s):=\lambda T J_1+\bigl(s+\beta H(cT\xi)\bigr)J_2=
\left\{
\begin{array}{ll}
\displaystyle \lambda T J_1+sJ_2=:A_T\equiv A_T(s),\ \ &\xi<0,\\[.3pc]
\displaystyle \lambda T J_1+\tilde{s}J_2=:\tilde{A}_T\equiv
A_T(\tilde{s}),\ \ &\xi>0,
\end{array}\right.
\end{equation}
\begin{equation}\label{eq:J12}
J_1:=\left(
\begin{array}{cc}
1& -1\\[.3pc]
1& -1
\end{array}
\right),\qquad J_2:=\left(
\begin{array}{cr}
1& 0\\[.3pc]
0& -1
\end{array}
\right).
\end{equation}
Note that
\begin{equation}\label{Jto1}
J_1\bfone={\bf 0}\myp,\qquad
J_1\tilde{\bfone}=2\cdot\mynn\bfone\myp,\qquad
J_2\myp\bfone=\tilde{\bfone}, \qquad
J_2\myp\tilde{\bfone}=\bfone\myp,
\end{equation}
where ${\bf 0}:=(0,0)'$. Hence (see (\ref{eq:A0}))
\begin{equation}\label{Ato1}
A_T(s)\mypp\bfone=s\mypp\tilde{\bfone},\qquad
A_T(s)\mypp\tilde{\bfone}=(s+2\lambda T)\mypp\bfone\myp.
\end{equation}

Let us set
\begin{equation}\label{eq:kappa}
\kappa\equiv\kappa(s):=\sqrt{s\mypp(s+2\lambda T)}\,,\qquad
\tilde\kappa:=\kappa(\tilde{s})=\sqrt{\tilde{s}\mypp(\tilde{s}+2\lambda
T)}\,.
\end{equation}
Using formulas (\ref{Ato1}) and (\ref{eq:kappa}), it is easy to
check that the matrix $A_T(s)$ has the eigenvalues $\pm\kappa(s)$
with the corresponding eigenvectors
\begin{equation}\label{eq:eigen-}
\bfa_\pm\equiv
\bfa_\pm(s):=\pm\kappa\myp\bfone+s\myp\tilde{\bfone}\myp,\qquad
A_T\myp\bfa_\pm=\pm\kappa\myp\bfa_\pm\,.
\end{equation}
In particular, relations (\ref{eq:eigen-}) imply that the
exponential of $A_T(s)$ can be represented as follows:
\begin{equation}\label{expA}
\rme^{A_T\xi}
={\textstyle\frac{1}{2}}\,\rme^{\kappa\xi}\bigl(I+\kappa^{-1}\myn
A_T\bigr)
+{\textstyle\frac{1}{2}}\,\rme^{-\kappa\xi}\bigl(I-\kappa^{-1}\myn
A_T\bigr),
\end{equation}
where $I$ is the identity matrix.

Recall that we are looking for a solution to the boundary value
problem (\ref{sysw-m0})\myp--\myp(\ref{limw-1}) continuous at the
origin. The following lemma gives an explicit form of such a
solution.
\begin{lemma}\label{lm:kappa}
For each $s>0$\myp\textup{,} the differential equation
\textup{(\ref{sysw-m0})} subject to the boundary conditions
\textup{(\ref{limw-1})} has the unique continuous solution given by
\begin{equation}\label{solw}
\bfw_T(\xi,s)=\left\{\!
\begin{array}{rl}
\displaystyle - {\rm e}^{\kappa \xi}\frac{\beta
s^{-1}}{s\tilde{\kappa}+\tilde{s}\kappa}\,(\kappa\myp
\bfone+s\myp\tilde{\bfone})+ s^{-1}\bfone\myp,\quad
&\xi\le0\myp,\\[.9pc]
\displaystyle \rme^{-\tilde{\kappa} \xi} \frac{\beta\mypp
\tilde{s}^{\myp-1}}{s\tilde{\kappa}+\tilde{s}\kappa}\,
(\tilde{\kappa}\myp\bfone-\tilde{s}\mypp\tilde{\bfone})+
\tilde{s}^{\myp-1}\bfone\myp,\quad & \xi\ge0\myp.
\end{array}
\right.
\end{equation}
In particular\textup{,}
\begin{equation}\label{w+00}
\bfw_T(0,s)=\frac{(\tilde{\kappa}+\kappa)\myp\bfone
+(s-\tilde{s})\myp\tilde{\bfone}}{s\tilde{\kappa}+\tilde{s}\kappa}\,.
\end{equation}
\end{lemma}

\proof Observe that the step function
$\bfw_T^*(\xi,s):=\bigl(s+\beta H(cT\xi)\bigr)^{-1}\bfone$ is a
particular solution of the equation (\ref{sysw-m0}) for each $s>0$
and all $\xi\ne0$. Indeed, the function $\bfw_T^*(\cdot,s)$ is
piecewise constant outside zero, hence
$(\partial/\partial\xi)\myp\bfw_T^*(\xi,s)=0$ ($\xi\ne0$), whereas,
due to (\ref{eq:A0}) and (\ref{Jto1}),
\[
\calA_T(\xi,s)\mypp\bfw_T^*(\xi,s)=\lambda T \bigl(s+\beta
H(cT\xi)\bigr)^{-1} J_1\bfone+J_2\myp \bfone\equiv \tilde{\bfone}.
\]

Therefore, a general solution of the linear differential equation
(\ref{sysw-m0}) can be represented in the form (see (\ref{eq:A0}))
\begin{equation}\label{eq:bfw}
\bfw_T(\xi,s)=\left\{\begin{array}{ll}
\rme^{A_T\xi}\mypp \bfc(s)+s^{-1}\bfone\myp,\quad &\xi<0,\\[.3pc]
\rme^{\tilde{A}_T\xi}\,
\tildebfc(s)+\tilde{s}^{\myp-1}\bfone\myp,\quad &\xi>0,
\end{array}
\right.
\end{equation}
with arbitrary vectors $\bfc(s)$, $\tilde\bfc(s)$ (which may also
depend on $T$). Let us now find suitable $\bfc(s)$ and
$\tildebfc(s)$ so that the solution $\bfw_T(\cdot,s)$ would satisfy
the required boundary conditions at infinity and the continuity
condition at zero. From the representation (\ref{eq:bfw}) it is
clear that conditions (\ref{limw-1}) are satisfied if and only if
\begin{equation}\label{eq:exp-lim}
\lim_{\xi\to-\infty} \rme^{A_T\xi}\mypp \bfc(s)={\bf 0}\mypp,\qquad
\lim_{\xi\to+\infty} \rme^{\tilde{A}_T\xi}\mypp \tildebfc(s)={\bf
0}\mypp.
\end{equation}
Recalling that $\tilde{A}_T(s)=A_T(\tilde{s})$ and using the
exponential formula (\ref{expA}), it is easy to see that conditions
(\ref{eq:exp-lim}) are reduced to the equations
\begin{equation*}
(I-\kappa^{-1}\myn A_T)\,\bfc(s)={\bf 0}\mypp,\qquad
(I+\tilde{\kappa}^{\myp-1}\tilde{A}_T)\,\tildebfc(s)={\bf 0}\mypp,
\end{equation*}
which implies that $\bfc(s)$ and $\tildebfc(s)$ are eigenvectors of
the matrices $A_T$ and $\tilde{A}_T$, respectively, with the
corresponding eigenvalues $\kappa$ and $-\tilde{\kappa}$. On account
of formulas (\ref{eq:eigen-}), this immediately gives
$\bfc(s)=C(s)\,\bfa_+$\myp,
\,$\tildebfc(s)=\tilde{C}(s)\,\tilde{\bfa}_-$\myp, with some
real-valued functions $C(s)$, $\tilde{C}(s)$. Therefore, after the
substitution of expressions (\ref{eq:eigen-}), formula
(\ref{eq:bfw}) takes the form
\begin{equation}\label{eq:bfw1}
\bfw_T(\xi,s)=\left\{\begin{array}{rl} \rme^{\kappa\xi}\mypp
C(s)\myp (\kappa\myp\bfone+s\myp\tilde{\bfone}) +s^{-1}\bfone\myp,
\ \ \quad &\xi< 0\myp,\\[.3pc]
-\rme^{-\tilde{\kappa}\xi}\mypp\tilde{C}(s)\myp
(\tilde\kappa\myp\bfone-\tilde{s}\myp\tilde{\bfone})+\tilde{s}^{\myp-1}\bfone\myp,\
\ \quad &\xi> 0\myp.
\end{array}
\right.
\end{equation}
Furthermore, taking into account the continuity of $\bfw_T(\cdot,
s)$ at zero, from (\ref{eq:bfw1}) we
have
$$
C(s)\mypp (\kappa\myp\bfone+s\myp\tilde{\bfone})
+s^{-1}\bfone=\tilde{C}(s)\mypp
(-\tilde\kappa\myp\bfone+\tilde{s}\myp\tilde{\bfone})+\tilde{s}^{\myp-1}\bfone\myp,
$$
whence, by equating the coefficients of $\bfone$ and
$\tilde{\bfone}$ on the left- and right-hand sides, we obtain
\begin{equation*}
\left\{
\begin{aligned}
&C(s)\mypp\kappa+s^{-1}=-\tilde{C}(s)\mypp\tilde{\kappa}+\tilde{s}^{\myp-1},\\
&C(s)\myp s=\tilde{C}(s)\mypp \tilde{s}\,.
\end{aligned}
\right.
\end{equation*}
Solving this system of equations we find
\begin{equation*}
C(s)=\frac{-\beta s^{-1}} {s\tilde{\kappa}+\tilde{s}\kappa}\,,\qquad
\tilde{C}(s)=\frac{-\beta\myp
\tilde{s}^{\myp-1}}{s\tilde{\kappa}+\tilde{s}\kappa}\,,
\end{equation*}
and the substitution of these expression into (\ref{eq:bfw1}) yields
the required formula (\ref{solw}).

Finally, the expression (\ref{w+00}) for $\bfw_T(0,s)$ follows from
(\ref{solw}) by setting $\xi=0$ and using that $\beta=\tilde{s}-s$
(see Lemma \ref{lm:tilde}). This completes the proof of Lemma
\ref{lm:kappa}.
\endproof

\begin{lemma}\label{lm:w-pm(0)}
The components $w_T^\pm(0,s)$ \textup{(}see
\textup{(\ref{w+00})}\textup{)} are explicitly given by the
expressions
\begin{align}\label{w+0}
w_T^+(0,s)&=\frac{2}{\tilde{\kappa}+\tilde{s}}+ \frac{2\lambda
T}{(\kappa+s)(\tilde{\kappa}+\tilde{s})}\,,\\
\label{w-0} w_T^-(0,s)&=\frac{2}{\kappa+s}+\frac{2\lambda
T}{(\kappa+s)(\tilde\kappa+\tilde s)}\,.
\end{align}
\end{lemma}

\proof From the vector expression (\ref{w+00}) we have
\begin{equation}\label{w+00-}
w_T^\pm(0, s)=\frac{\tilde{\kappa}\mp\tilde{s} +\kappa\pm
s}{s\tilde{\kappa}+\tilde{s}\kappa}\,.
\end{equation}
Note that, according to (\ref{eq:kappa}),
\begin{equation}\label{w+00--}
\kappa^2-s^2=2\lambda Ts,\qquad \tilde\kappa^2-\tilde
s^{\myp2}=2\lambda T\tilde s,
\end{equation}
hence the expression (\ref{w+00-}) may be rewritten as
\begin{align}
\notag w_T^\pm(0,s)&=\frac{1}{s\tilde\kappa+\tilde
s\kappa}\left(\frac{\kappa^2-s^2}{\kappa\mp
s}+\frac{\tilde\kappa^2-\tilde s^2}{\tilde\kappa\pm\tilde
s}\right)\\
\label{kappapms}
&=\frac{2\lambda T}{s\tilde\kappa+\tilde
s\kappa}\left(\frac{s}{\kappa\mp s}+\frac{\tilde
s}{\tilde\kappa\pm\tilde s}\right) =\frac{2\lambda T}{(\kappa\mp
s)(\tilde\kappa\pm \tilde{s})}\,,
\end{align}
which is equivalent to (\ref{w+0}), (\ref{w-0}); for instance, for
$w_T^+(0,s)$ (corresponding to the choice of the upper sign in $\pm$
and $\mp$), from formula (\ref{kappapms}) we obtain, again using
(\ref{w+00--}),
\begin{align*}
w_T^+(0,s)&=\frac{2\lambda
T}{(\kappa-s)(\tilde\kappa+\tilde{s})}=\frac{2\lambda
T(\kappa+s)}{(\kappa^2-s^2)(\tilde\kappa+\tilde{s})}
=\frac{\kappa+s}{s\myp(\tilde\kappa+\tilde{s})}\\
&=\frac{2}{\tilde\kappa+\tilde{s}}
+\frac{\kappa-s}{s\myp(\tilde\kappa+\tilde{s})}
=\frac{2}{\tilde\kappa+\tilde{s}}+\frac{2\lambda
T}{(\kappa+s)(\tilde\kappa+\tilde{s})}\,,
\end{align*}
in agreement with (\ref{w+0}). Thus, Lemma \ref{lm:w-pm(0)} is
proved.
\endproof

\begin{lemma}\label{lem33}
Let the function $\varphi_T(t)$ be defined by
\textup{(\ref{eq:phi_T})}. Then\myp\textup{,} for each
$s>0$\myp\textup{,}
\begin{equation}\label{A}
\int_0^\infty \rme^{-s\myp t}\mypp\varphi_T(t)\,\rmd
t=\frac{1}{\kappa+s}\,,\qquad \int_0^\infty \rme^{-s\myp t}\mypp
\rme^{-\beta t} \varphi_T(t)\,\rmd
t=\frac{1}{\tilde{\kappa}+\tilde{s}}\,,
\end{equation}
and
\begin{equation}\label{B}
\int_0^\infty \rme^{-s\myp t}\left(\int_0^t\rme^{-\beta
y}\mypp\varphi_T(y)\mypp\varphi_T(t-y)\,\rmd y\right)\rmd
t=\frac{1}{(\kappa+s)(\tilde{\kappa}+\tilde{s})},
\end{equation}
where $\tilde{s}=s+\beta$ and $\kappa=\kappa(s)$\textup{,}
$\tilde\kappa=\kappa(\tilde s)$ are defined in
\textup{(\ref{eq:kappa})}.
\end{lemma}

\proof Inserting (\ref{eq:phi_T}) and changing the order of
integration, we obtain
\begin{align}
\notag
\int_0^\infty \rme^{-s\myp t}\mypp\varphi_T(t)\,\rmd t
&=\frac{1}{4\pi\lambda T}\int_0^\infty\frac{1-\rme^{-2\lambda
Tu}}{u^{3/2}}\left(\int_u^\infty\frac{\rme^{-s\myp
t}}{\sqrt{t-u}}\,\rmd t\right)\rmd u\\
\label{eq:-3/2} &=\frac{\Gamma(\mbox{$\frac12$})}{4\pi\lambda
T\sqrt{s}}\int_0^\infty \rme^{-s\myp u}\,\bigl(1-\rme^{-2\lambda
Tu}\bigr)\,u^{-3/2}\,\rmd{u}.
\end{align}
Integration by parts via $u^{-3/2}\,\rmd{u}=-2\,\rmd(u^{-1/2})$
yields the right-hand side of (\ref{eq:-3/2}) in the form
\begin{align*}
\frac{1}{2\lambda T\sqrt{\pi s}}\int_0^\infty
u^{-1/2}\bigl(\rme^{-(s+2\lambda T)\myp u}-\rme^{-s\myp
u}\bigr)\,\rmd{u}&=\frac{\bigl(\sqrt{s+2\lambda
T}-\sqrt{s}\,\bigr)\,\Gamma(\mbox{$\frac12$})}{2\lambda T\sqrt{\pi
s}} =\frac{1}{\kappa+s}\,,
\end{align*}
and the first formula in (\ref{A}) is proved. The second one readily
follows by the shift $\tilde{s}=s+\beta$.

Furthermore, using the convolution property of the Laplace
transform, the left-hand side of (\ref{B}) is reduced to the product
\begin{align*}
\int_0^\infty &\rme^{-s\myp t}\,\rme^{-\beta
t}\mypp\varphi_T(t)\,\rmd{t}\,\times \int_0^\infty
\rme^{-st}\mypp\varphi_T(t)\,\rmd t
=\frac{1}{(\tilde{\kappa}+\tilde{s})(\kappa+s)}\,,
\end{align*}
according to (\ref{A}), which completes the proof of the lemma.
\endproof

Combining Lemmas \ref{lm:w-pm(0)} and \ref{lem33} and using the
uniqueness theorem for the Laplace transform (\ref{eq:w}), we obtain
\begin{align*}
v_T^{\pm}(0,t)= \bigl(1+\rme^{-\beta t}\mp 1\pm \rme^{-\beta
t}\bigr)\mypp \varphi_T(t)+2\lambda T\int_0^t\rme^{-\beta
y}\myp\varphi_T(y)\mypp\varphi_T(t-y)\,\rmd{y}\myp.
\end{align*}
In particular, setting $t=1$ (see (\ref{eq:v})) and recalling the
definition (\ref{eq:psi_T}) of the function $\psi_T$, we get
\begin{equation*}
\EE\bigl[\rme^{-\beta\eta^{\pm}_T(0)}\bigr]=
\bigl(1+\rme^{-\beta}\mp 1\pm \rme^{-\beta}\bigr)\mypp \varphi_T(1)+
\int_0^1 \rme^{-\beta y}\myp\psi_T(y)\,\rmd{y}\myp,
\end{equation*}
and it is evident (in view of the uniqueness theorem for
Laplace transform)
that the distribution of $\eta^\pm_T(0)$ is given by formula
(\ref{eta0a}).

Finally, the result (\ref{eta00-}) for $\eta_T(0)$ readily follows
from (\ref{eta0a}) and the decomposition
\begin{equation}\label{eq:1/2+1/2}
\PP\bigl\{\eta_T(0)\in
{\rmd}y\bigr\}=\frac12\,\PP\bigl\{\eta^{+}_T(0)\in
{\rmd}y\bigr\}+\frac12\,\PP\bigl\{\eta^{-}_T(0)\in
{\rmd}y\bigr\}\qquad(0\le y\le 1).
\end{equation}
Thus, the proof of Theorem \ref{th1} is completed.

\section{Proof of Theorem \ref{th2}}\label{sec:th2}

The plan of the proof below is to calculate the Laplace transform
(see (\ref{eq:v}) and (\ref{eq:w}))
\begin{equation}\label{eq:w-new}
w_T^\pm(\xi,s)=\int_0^\infty\rme^{-st}\, \EE\bigl[\rme^{-\beta
t\mypp\eta^\pm_{Tt}(cT\xi)}\bigr]\,\rmd{t}
\end{equation}
from the explicit (hypothetical) distribution of $\eta_T^{\pm}(x)$
given by formula (\ref{eq:Y-}), and to verify that the result
coincides with formulas (\ref{solw}) obtained in Lemma \ref{solw}.
The claim of Theorem \ref{th2} will then follow by the uniqueness
theorem for Laplace transform. To be specific, we will focus on the
$w_T^{+}$ case, the proof for $w_T^{-}$ being similar.

\begin{remark}
In Appendix \ref{sec:A2} we will give an alternative proof based on
probabilistic arguments, by reducing the general case
$\eta_T^{\pm}(x)$ to $\eta_T^{\pm}(0)$ via conditioning on the
hitting time of the origin. That proof will explain how formulas
(\ref{eq:Y-}) can be derived (rather than verified); however, in so
doing the prior knowledge of the distribution of $\eta_T^{\pm}(0)$
(provided by Theorem \ref{th1}) is essential.
\end{remark}

Due to the space-time change $(x,T)\mapsto(cT\xi,Tt)$ used in
(\ref{eq:w-new}), the time threshold $T_0=|x|/c$ becomes
$T_0=T|\xi|$, whereas the former condition $T>T_0$ is converted into
$t>|\xi|$. As a first step in the proof, using the probability
distribution proposed by the theorem (including its part (a)) we can
represent the Laplace transform of $t\myp\eta_{Tt}^{+}(cT\xi)$ as
\begin{equation}\label{eq:LT-J}
\EE\bigl[\rme^{-\beta t\myp\eta_{Tt}^{+}(cT\xi)}\bigr]
=\left\{\begin{array}{cl} \displaystyle 1,&t\le|\xi|,\\
\displaystyle\sum_{i=1}^5 \calJ_T^{(i)}(\xi,t),\ \ &t>|\xi|,
\end{array} \right.
\end{equation}
where $\calJ_T^{(i)}(\xi,t)$ ($i=1,\dots,5$) arise from the three
parts on the right-hand side of the representation (\ref{eq:Y-}),
with the last two further subdivided each into two terms, according
to (\ref{eq:mu}) and (\ref{eq:varPsi}). More precisely, using the
scaling property (\ref{eq:phi-scaling}) of the function $\varphi_T$
and making the substitutions $y\mapsto ty$ and $u\mapsto T(u+|\xi|)$
wherever appropriate, the functions $\calJ_T^{(i)}(\xi,t)$ can be
expressed as
\begin{gather}\label{eq:J-1}
\calJ_T^{(1)}(\xi,t):= T\int_{t}^\infty
\!Q^{+}_{cT|\xi|}(Tu)\,\rmd{u},\\[.3pc]
\label{eq:J-2} \calJ_T^{(2)}(\xi,t+|\xi|):= 2\mypp\rme^{-\beta
 t-\lambda T|\xi|}\, \varphi_{T}(t),\\[.3pc]
\label{eq:J-3} \calJ_T^{(3)}(\xi,t+|\xi|):=\rme^{-\lambda
T|\xi|}\int_0^{t}\rme^{-\beta y}
\,\psi_{Tt}\!\left(\frac{y}{t}\right)\frac{\rmd{y}}{t},\\[.3pc]
\label{eq:J-4} \calJ_T^{(4)}(\xi,t+|\xi|):=2\myp T\int_0^{t}
\rme^{-\beta y}\mypp\varphi_{T}(y)\,
Q^{+}_{cT|\xi|}\bigl(T(t+|\xi|-y)\bigr)\,\rmd{y},\\[.3pc]
\label{eq:J-5} \calJ_T^{(5)}(\xi,t+|\xi|):=T\int_0^{t} \rme^{-\beta
y}\left(\int_{0}^{t-y}
\!Q^{+}_{cT|\xi|}\bigl(T(u+|\xi|)\bigr)\,\psi_{T(t-u)}
\!\left(\frac{y}{t-u}\right)\frac{\rmd{u}}{t-u}\right)\rmd{y}.
\end{gather}
Consequently, from (\ref{eq:w-new}) and (\ref{eq:LT-J}) we get
\begin{equation}\label{eq:ww}
w^{+}_T(\xi,s)=\frac{1}{s}\,\bigl(1-\rme^{-|\xi| s}\bigr)
+\sum_{i=1}^5\int_{|\xi|}^\infty\rme^{-st} \myp\calJ_T^{(i)}(\xi,t)\
\rmd{t}\myp.
\end{equation}

Let us now calculate the Laplace transform (with respect to $t$) of
each of the terms $\calJ_T^{(i)}(\xi,t)$ ($i=1,\dots,5$). In so
doing, the next formula will be useful,
\begin{equation}\label{eq:Q+Laplace}
T\int_{|\xi|}^\infty \rme^{-st} \mypp
Q_{cT|\xi|}^{+}(Tt)\,\rmd{t}=\rme^{\xi\kappa}-\rme^{(s+\lambda
T)\myp\xi}\qquad (\xi<0),
\end{equation}
where $\kappa=\sqrt{s(s+2\lambda T)}$ \,(see (\ref{eq:kappa})),
which immediately follows from the definition (\ref{eq:Q+})
according to \cite[29.3.96, p.\,1027]{AS}.

\begin{remark}
An analogous formula for $Q^{-}$ (needed for the proof in the case
of $w^{-}$) follows from (\ref{eq:Q-}) by applying \cite[29.3.93,
29.3.96, p.\,1027]{AS}).
\end{remark}

(i) \,From (\ref{eq:J-1}) we obtain, integrating by parts and using
(\ref{eq:Q+Laplace}),
\begin{align}
\notag \int_{|\xi|}^\infty
\rme^{-st}\myp\calJ_T^{(1)}(\xi,t)\,\rmd{t} &=Ts^{-1}\myp
\rme^{-s|\xi|}\int_{|\xi|}^\infty
Q^{+}_{cT|\xi|}(Tt)\,\rmd{t}-Ts^{-1}\!\int_{|\xi|}^\infty
\rme^{-st}\mypp Q^{+}_{cT|\xi|}(Tt)\,\rmd{t}\\[.3pc]
&=\frac{1}{s}\mypp\bigl(\rme^{s\xi}-\rme^{(s+\lambda
T)\mypp\xi}\bigr)
-\frac{1}{s}\mypp\bigl(\rme^{\kappa\xi}-\rme^{(s+\lambda
T)\mypp\xi}\bigr)=\frac{1}{s}\mypp\bigl(\rme^{s\xi}-\rme^{\kappa\xi}\bigr).
\label{eq:J1=}
\end{align}

\smallskip
(ii) \,After the substitution $t\mapsto t+|\xi|$, from
(\ref{eq:J-2}) we get, using formula (\ref{A}) in Lemma \ref{lem33},
\begin{align}
\int_{|\xi|}^\infty \rme^{-st}\myp\calJ_T^{(2)}(\xi,t)\,\rmd{t}
&=2\mypp\rme^{(s+\lambda T)\myp \xi} \int_{0}^\infty
\rme^{-(s+\beta)\myp t}\myp \varphi_{T}(t)\,\rmd{t}
=2\mypp\rme^{(s+\lambda T)\myp \xi}
\mypp\frac{1}{\tilde{\kappa}+\tilde{s}}\,.\label{eq:J2=}
\end{align}

\smallskip
(iii) \,Likewise, from (\ref{eq:J-3}) we obtain, recalling the
definition (\ref{eq:psi_T}) of the function $\psi_T$ and again using
the scaling property (\ref{eq:phi-scaling}),
\begin{align}
\notag \int_{|\xi|}^\infty
\rme^{-st}\myp\calJ_T^{(3)}(\xi,t)\,\rmd{t} &=\rme^{(s+\lambda
T)\myp \xi}\int_{0}^\infty \rme^{-st}\left(\int_0^{t}\rme^{-\beta y}
\,\psi_{Tt}\!\left(\frac{y}{t}\right)\frac{\rmd{y}}{t}\right)\rmd{t}\\
\notag &=\rme^{(s+\lambda T)\myp \xi}\,2\lambda T\int_{0}^\infty
\rme^{-st}\left(\int_0^{t}\rme^{-\beta y}
\,\varphi_{T}(y)\,\varphi_T(t-y)\,\rmd{y}\right)\rmd{t}\\
\label{eq:J3=}  &=\rme^{(s+\lambda T)\myp \xi} \,\frac{2\lambda
T}{(\kappa+s)(\tilde{\kappa}+\tilde{s})}\,,
\end{align}
as follows from formula (\ref{B}) in Lemma \ref{lem33}.

\smallskip
(iv) \,Similarly, taking advantage of the convolution theorem, the
Laplace transform of (\ref{eq:J-4}) can be written as
\begin{align}
\notag \int_{|\xi|}^\infty \rme^{-st}\,\calJ_T^{(4)}(\xi,t)\,\rmd{t}
&=\int_{0}^\infty \rme^{-st}\left(2\myp T\int_0^{t} \rme^{-\beta
y}\mypp\varphi_{T}(y)\,
Q^{+}_{cT|\xi|}\bigl(T(t+|\xi|-y)\bigr)\,\rmd{y}\right)\rmd{t}\\
\notag&=2\int_{0}^\infty \rme^{-st}\mypp\rme^{-\beta
t}\mypp\varphi_{T}(t)\,\rmd{t}\,\times \rme^{s\myp\xi}\mypp
T\int_{0}^\infty
\rme^{-st}\mypp Q^{+}_{cT|\xi|}\bigl((t+|\xi|)T\bigr)\,\rmd{t}\\
\label{eq:J4=}
&=\frac{2}{\tilde{\kappa}+\tilde{s}}\mypp\bigl(\rme^{\kappa\xi}-\rme^{(s+\lambda
T)\myp\xi}\bigr),
\end{align}
according to formulas (\ref{A}) and (\ref{eq:Q+Laplace}).

\smallskip
(v) \,Interchanging the integrations, we can rewrite (\ref{eq:J-5})
in the form
\begin{align*}
\calJ_T^{(5)}(\xi,t+|\xi|) &= T\int_{0}^{t}
Q^{+}_{cT|\xi|}\bigl(T(u+|\xi|)\bigr)\left(\int_0^{t-u}\rme^{-\beta
y}\mypp\psi_{T(t-u)}\left(\frac{y}{t-u}\right)\,\frac{\rmd{y}}{t-u}\right)\rmd{u},
\end{align*}
hence, by the convolution theorem, the Laplace transform of
$\calJ_T^{(5)}(\xi,t)$ is reduced to
\begin{align}
\notag \int_{|\xi|}^\infty \rme^{-st}\,\calJ_T^{(5)}(\xi,t)\,\rmd{t}
&= \rme^{s\myp\xi}\,T\int_{0}^\infty \rme^{-st}\mypp
Q^{+}_{cT|\xi|}\bigl(T(t+|\xi|)\bigr)\,\rmd{t}\,\times
\int_{0}^\infty \rme^{-st}\left(\int_0^t \rme^{-\beta y}\mypp
\psi_{Tt}\left(\frac{y}{t}\right)\,\frac{\rmd{y}}{t}\right)\rmd{t}\\[.3pc]
\label{eq:J5=} &=\bigl(\rme^{\kappa\xi}-\rme^{(s+\lambda
T)\myp\xi}\bigr)\frac{2\lambda
T}{(\kappa+s)(\tilde{\kappa}+\tilde{s})}\,,
\end{align}
as was shown in (\ref{eq:J3=}) and (\ref{eq:J4=}).

\smallskip
Finally, substituting the results (\ref{eq:J1=}), (\ref{eq:J2=}),
(\ref{eq:J3=}), (\ref{eq:J4=}) and (\ref{eq:J5=}) into formula
(\ref{eq:ww}) and recalling the expressions (\ref{w+0}), (\ref{w-0})
for $w_T^{\pm}(0,s)$, we get
\begin{align*}
w_T^{+}(\xi,s)&=\frac{1}{s}\,\bigl(1-\rme^{\kappa\xi}\bigr)
+\rme^{\kappa\xi}\left(\frac{2}{\tilde{\kappa}+\tilde{s}}+\frac{2\lambda
T}{(\kappa+s)(\tilde{\kappa}+\tilde{s})}\right)\\[.2pc]
&=\frac{1}{s}\,\bigl(1-\rme^{\kappa\xi}\bigr)+\rme^{\kappa\xi}\,w_T^{+}(0,s)=
\rme^{\kappa\xi}\mynn\left(w_T^{+}(0,s)-\frac{1}{s}\right)+\frac{1}{s}\,,
\end{align*}
which is consistent with the expression (\ref{solw}) for
$w_T^{+}(\xi,s)$ obtained in Lemma \ref{lm:kappa}. Thus, the proof
of Theorem \ref{th2} is completed.

\section{Proof of Theorem \ref{th3}}\label{sec:th3}

It suffices to prove the theorem for the conditional versions
$\eta_T^\pm(x)$ only; indeed, since the latter have the same
distributional limit, the result for $\eta_T(x)$ will readily follow
(cf.\ (\ref{eq:1/2+1/2})).

In the next lemma, we find the Laplace transform for a suitable
parametric family $Y_a(t)$ extending the random variables $Y_a$
introduced in Section~\textup{\ref{sec:Results}} (see (\ref{eq:Y0}),
(\ref{eq:m_a}) and (\ref{eq:f_a})). Recall that $\tilde s=s+\beta$.

\begin{lemma}\label{lm:theta}
For any $a\ge0$ and $t>0$\myp\textup{,} set \,$Y_a(t):=t\,
Y_{a/\sqrt{t}}$\mypp. Then\myp\textup{,} for any $s>0$ and
$\beta>0$\myp\textup{,} we have
\begin{align}\label{eq:l1pm}
\int_0^\infty \rme^{-s\myp t}\, \EE \bigl[\rme^{-\beta \myp
Y_a(t)}\bigr]\,\rmd t&
=\rme^{-a\sqrt{2s}}\left(\frac{1}{\sqrt{s\tilde{s}}}-\frac{1}{s}\right)+\frac{1}{s}\,,\\
\label{eq:l2pm} \int_0^\infty \rme^{-s\myp t}\, \EE
\bigl[\rme^{-\beta \myp (t-Y_a(t))}\bigr]\,\rmd t &
=\rme^{-a\sqrt{2\tilde{s}}}\left(\frac{1}{\sqrt{s\tilde{s}}}-\frac{1}{\tilde{s}}\right)+
\frac{1}{\tilde{s}}\,.
\end{align}
In particular\textup{,} for $a=0$
\begin{equation}\label{eq:l1}
\int_0^\infty \rme^{-s\myp t}\, \EE \bigl[\rme^{-\beta \myp
Y_0(t)}\bigr]\,\rmd{t} =
\frac{1}{\sqrt{s\tilde
s}}\,.
\end{equation}
\end{lemma}

\proof It is sufficient to prove formula (\ref{eq:l1pm}) only;
indeed,
\begin{equation}\label{eq:-beta}
\int_0^\infty \rme^{-s\myp t}\, \EE\myp \bigl[\rme^{-\beta
(t-Y_a(t))}\bigr]\mypp\rmd t =\int_0^\infty \rme^{-\tilde{s}\myp
t}\,\EE \bigl[\rme^{\beta\myp Y_a(t)}\bigr]\mypp\rmd{t},
\end{equation}
hence the left-hand side of (\ref{eq:l2pm}) can be computed using
(\ref{eq:l1pm}) by changing $s$ to $\tilde{s}$ and $\beta$ to
$-\beta$, which amounts to interchanging the symbols $s$ and
$\tilde{s}$ in (\ref{eq:l1pm}), thus leading to formula
(\ref{eq:l2pm}). (Note that the right-hand side of (\ref{eq:-beta})
is well defined since $Y_a(t)\le t$ and so $\rme^{-(s+\beta)\myp
t}\,\EE\myp[\rme^{\beta \myp Y_a(t)}]\le \rme^{-st}$.)

Now, if $a=0$ then $Y_0(t)=t\mypp Y_0$, where $Y_0$ has the arcsine
distribution with the density (\ref{eq:arcsine}), hence the
left-hand side of (\ref{eq:l1}) is reduced to
\begin{equation}\label{eq:conv}
\int_0^\infty \rme^{-s\myp t}\left(\frac{1}{\pi}\int_0^t
\frac{\rme^{-\beta y}}{\sqrt{y(t-y)}}\,\rmd{y}\right)\rmd{t}\myp.
\end{equation}
The internal integral here can be interpreted as the convolution
$(f_1* f_2)(t)$ of the functions $f_1(t)=\rme^{-\beta t} \mypp
t^{-1/2}$ and $f_2(t)= t^{-1/2}$, hence the Laplace transform
(\ref{eq:conv}) reduces to the product
$$
\frac{1}{\pi}\int_0^\infty \rme^{-s\myp t}\,\rme^{-\beta t}
\,t^{-1/2}\,\rmd t \int_0^\infty \rme^{-s\myp t}\,t^{-1/2}\,\rmd t
=\frac{\Gamma(\mbox{$\frac12$})}{\pi\sqrt{s+\beta}}\cdot
\frac{\Gamma(\mbox{$\frac12$})}{\sqrt{s}}
=\frac{1}{\sqrt{\tilde{s}\myp s}}\,,
$$
and the required formula (\ref{eq:l1}) follows.

If $a>0$ then, noting that $q_{a/\sqrt{t}}(u)=t\mypp q_a(ut)$ and
using (\ref{eq:Y0}), (\ref{eq:m_a}) and (\ref{eq:f_a}), we have
\begin{equation}\label{eq:LT-Y}
\EE \bigl[\rme^{-\beta\myp Y_a(t)}\bigr]=\int_t^\infty
q_{a}(u)\,\rmd{u}+\int_0^t \rme^{-\beta\myp y} \left(\int_0^{t-y}
\frac{q_{a}(u)}{t-u}\,p_{\mathrm{as}}\left(\frac{y}{t-u}\right)\rmd{u}\right)
\rmd{y}.
\end{equation}
Interchanging the order of integration and making the substitution
$y=z(t-u)$, we can rewrite the second (iterated integral) term on
the right-hand side of (\ref{eq:LT-Y}) as
\begin{equation*}
\int_0^t q_{a}(u)\left(\int_0^1 \rme^{-\beta (t-u)\myp z}\,
p_{\mathrm{as}}(z)\,\rmd{z}\right)\rmd{u},
\end{equation*}
which can be viewed as the convolution $(q_{a}*\hat{p}_{\beta})(t)$,
where
\begin{equation}\label{eq:p-hat}
\hat{p}_{\beta}(t):=\int_0^1 \rme^{-\beta t\myp z}\,
p_{\mathrm{as}}(z)\,\rmd{z}=\EE \bigl[\rme^{-\beta \myp
Y_0(t)}\bigr].
\end{equation}
Returning to (\ref{eq:LT-Y}) and applying the Laplace transform
(with respect to the variable $t$), by the convolution theorem the
left-hand side of (\ref{eq:l1pm}) can be expressed as
\begin{equation}\label{eq:4.8}
\int_0^\infty \rme^{-s\myp t}\left(\int_t^\infty
q_{a}(u)\,\rmd{u}\right)\rmd{t} + \int_0^\infty \rme^{-s\myp t}\mypp
q_{a}(t)\,\rmd{t}\,\times \int_0^\infty \rme^{-s\myp t}\mypp
\hat{p}_{\beta}(t)\,\rmd{t}.
\end{equation}
Recall that, according to (\ref{eq:4.9}),
\begin{equation}\label{eq:q-hat}
\int_0^\infty \rme^{-s\myp t}\mypp
q_{a}(t)\,\rmd{t}=\rme^{-a\sqrt{2s}},
\end{equation}
whence, integrating by parts and using (\ref{eq:4.9}), we
obtain
\begin{equation}\label{eq:4.11} \int_0^\infty \rme^{-s\myp
t}\left(\int_t^\infty q_{a}(u)\,\rmd{u}\right)\rmd{t}
=\frac{1}{s}-\frac{1}{s}\,\rme^{-a\sqrt{2s}}.
\end{equation}
Furthermore, from (\ref{eq:p-hat}) and (\ref{eq:l1}) we have
\begin{equation}\label{eq:4.10}
\int_0^\infty \rme^{-s\myp t}\mypp
\hat{p}_{\beta}(t)\,\rmd{t}=\frac{1}{\sqrt{s\tilde{s}}}\,.
\end{equation}
As a result, substituting expressions (\ref{eq:q-hat}),
(\ref{eq:4.11}) and (\ref{eq:4.10}) into (\ref{eq:4.8}), we obtain
formula (\ref{eq:l1pm}).
\endproof

\proof[Proof of  Theorem \textup{\ref{th3}}] As $T\to\infty$, we
have $\xi:=(cT)^{-1}x=(\lambda T)^{-1/2}\bigl(a+o(1)\bigr)$, whereas
from (\ref{eq:kappa}) it follows that $\kappa(s)\sim (2\lambda T
s)^{1/2}$, $\tilde{\kappa}(s)\sim(2\lambda T \tilde{s})^{1/2}$.
Hence, from (\ref{solw}) we obtain, for $\xi\le 0$, $a\le 0$,
\begin{equation*}
\lim_{T\to\infty} w_T^\pm(\xi, s)= -
\rme^{a\sqrt{2s}}\frac{\beta\mypp
s^{-1}}{\sqrt{\tilde{s}}\,(\sqrt{\tilde{s}}+\sqrt{s}\,)}+\frac{1}{s}
=\rme^{a\sqrt{2s}}\left(\frac{1}{\sqrt{s\tilde{s}}}-\frac{1}{s}\right)+\frac{1}{s}\,,
\end{equation*}
and similarly, for $\xi\ge0$, $a\ge 0$,
\begin{equation*}
\lim_{T\to\infty} w_T^\pm(\xi, s)=
\rme^{-a\sqrt{2\tilde{s}}}\frac{\beta\,\tilde{s}^{\myp-1}}
{\sqrt{s}\,(\sqrt{s}+\sqrt{\tilde{s}}\,)}+\frac{1}{\tilde{s}}
=\rme^{-a\sqrt{2\tilde{s}}}\left(\frac{1}{\sqrt{s\tilde{s}}}-\frac{1}{\tilde{s}}\right)+
\frac{1}{\tilde{s}}\,.
\end{equation*}

Comparing these results with Lemma \ref{lm:theta}, by the continuity
theorem for Laplace transforms \cite{Feller} we conclude that, for
each $t>0$, the distribution of the random variable
$t\mypp\eta^\pm_{Tt}(x)$ (see (\ref{OT0}) and (\ref{eq:v}))
converges weakly, as $T\to\infty$, to the arcsine distribution
(\ref{eq:arcsine}) if $a=0$ and to the distribution of either
$Y_{-a}(t)$ if $a<0$ (see (\ref{eq:l1pm})) or $t-Y_a(t)$ if $a>0$
(see (\ref{eq:l2pm})). Specialized to the case $t=1$, this readily
gives the result of Theorem \ref{th3}. Thus the proof is completed.
\endproof

\section{Proof of Theorem \ref{th4}}\label{sec:th4}

Similarly to Section \ref{sec:th1}, let us set
$\bfw_T(\xi,s):=(w_T^{+}(\xi,s),w_T^{-}(\xi,s))'$ and rewrite
equations (\ref{sysw}), (\ref{wc*}) in the matrix form (cf.\
(\ref{sysw-m0}), (\ref{limw-1}))
\begin{gather}\label{sysw-m-}
\frac{\partial\myp\bfw_T(\xi,s;f) }{\partial \xi}= \calA_T(\xi,s;f)
\mypp \bfw_T(\xi,s;f)
-\tilde{\bfone},\qquad \xi\in\RR\setminus D_{\myn f}\myp,\\
\label{limw-} \lim_{\xi\to-\infty}\bfw_T(\xi,s;f)
=s^{-1}\myp\bfone\myp,\qquad \lim_{\xi\to+\infty}\bfw_T(\xi,s;f)
=\tilde{s}^{\myp-1}\myp\bfone\myp,
\end{gather}
with the matrix (cf.\ (\ref{eq:A0}))
\begin{equation*}
\calA_T(\xi,s;f):=\lambda T J_1+\bigl(s+\beta f(cT\xi)\bigr)\myp
J_2\myp,
\end{equation*}
where $J_1$ and $J_2$ are defined in (\ref{Jto1}). Let us set
\begin{equation}\label{w-w*}
\bfdelta_T(\xi,s):=\bfw_T(\xi,s;f)-\bfw_T(\xi, s;H),
\end{equation}
where $H$ is the Heaviside step function (cf.\ (\ref{eq:A0})). Owing
to the properties of the solution $\bfw_T(\xi,s;f)$ (see the end of
Section \ref{sec:FK}), the function $\bfdelta_T(\xi,s)$ is bounded
and continuous in $\xi\in\RR$ (for any fixed $T, \,s>0$). From the
relation (\ref{w-w*}) and equations (\ref{sysw-m-}) and
(\ref{limw-}), we obtain the differential equation
\begin{equation}\label{eq:du-*}
\frac{\partial\myp\bfdelta_T(\xi,s) }{\partial \xi}=
\calA_T(\xi,s;H)\mypp\bfdelta_T(\xi,s)+f_0(cT\xi)
\mypp\bbfw_T(\xi,s)\myp,\qquad \xi\in\RR\setminus (D_{\myn
f}\cup\{0\}),
\end{equation}
where $f_0:=f-H$ and $\bbfw_T:=\beta J_2\mypp\bfw_T$ (for short),
with the boundary conditions
\begin{equation}\label{limu*}
\lim_{\xi\to\pm\infty}\bfdelta_T(\xi,s)={\bf 0}\mypp.
\end{equation}
More explicitly, equation (\ref{eq:du-*}) splits into two equations
on the negative and positive half-lines:
\begin{align}\label{eqn:calB-minus*}
\frac{\partial \bfdelta_T(\xi,s)}{\partial
\xi}&=A_T\myp\bfdelta_T(\xi,s)+
f_0(cT\xi)\mypp \bbfw_T(\xi,s), \qquad \xi<0,\\
\label{eqn:calB-plus*} \frac{\partial
\bfdelta_T(\xi,s)}{\partial\xi} &=\tilde{A}_T\myp\bfdelta_T(\xi,s)+
f_0(cT\xi)\mypp \bbfw_T(\xi,s), \qquad \xi>0,
\end{align}
where $A_T\equiv A_T(s)=\lambda T J_1+s\myp J_2$, $\tilde{A}_T\equiv
A_T(\tilde{s})=\lambda T J_1+\tilde{s}\myp J_2$ (cf.\
(\ref{eq:A0})).

By the variation of constants, equation (\ref{eqn:calB-minus*}) is
equivalent to the integral equation
\begin{equation} \label{intsysu-}
\bfdelta_T(\xi,s)=\rme^{\xi A_T}\bfc_T+\int_0^\xi \rme^{(\xi-y)\myp
A_T}\myn f_0(cTy)\mypp\bbfw_T(y,s)\,\rmd{y},\qquad \xi\le0,
\end{equation}
where $\bfc_T\equiv\bfc_T(s)=\lim_{\xi\to0-} \bfdelta_T(\xi,s)$ is a
constant vector (for fixed $T$ and $s$). By the exponential formula
(\ref{expA}), equation (\ref{intsysu-}) takes the form
\begin{equation}\label{eq:D12*}
\bfdelta_T(\xi,s)=\frac{1}{2}\,\rme^{\kappa\xi}\Bigl[(I+\kappa^{-1}\myn
A_T)\mypp\bfc_T(s)+\bfq_T^{+}(\xi,s)\Bigr]
+\frac{1}{2}\,\rme^{-\kappa\xi}\Bigl[(I-\kappa^{-1}\myn A_T)\mypp
\bfc_T(s)+\bfq_T^{-}(\xi,s)\Bigr],
\end{equation}
where
\begin{equation}\label{eq:int_pm*}
\bfq_T^{\pm}(\xi,s):=(I\pm\kappa^{-1}\myn A_T)\int_0^\xi
\rme^{\mp\kappa y} f_0(cTy)\mypp\bbfw_T(y,s)\,\rmd{y},\qquad
\xi\le0.
\end{equation}

For fixed $s$ and $T$, we have $\bfq_T^{+}(\xi,s)
=\rme^{-\kappa\xi}\mypp o(1)$ as $\xi\to-\infty$. Indeed, via the
change of variables $z=y-\xi$ and applying Lebesgue's dominated
convergence theorem, we see that
\begin{align*}
\left|\int_0^\xi \rme^{-\kappa (y-\xi)}
f_0(cTy)\mypp\bbfw_T(y,s)\,\rmd{y}\right|= O(1)\int_0^\infty
\rme^{-\kappa z}\myp\bigl|f_0(cT(z+\xi))\bigr|\,\rmd{z}= o(1),\qquad
\xi\to-\infty,
\end{align*}
since $\bbfw_T$ and $f_0$ are bounded whereas $f_0(cT(z+\xi))\to0$
for each $z$, according to the hypothesis of Theorem \ref{th4}.
Hence, due to the boundary condition (\ref{limu*}) at $\xi=-\infty$,
equation (\ref{eq:D12*}) implies
\begin{equation}\label{eq:l'Hopital*}
\rme^{-\kappa\xi}\mypp\Bigl\{(I-\kappa^{-1}\myn A_T)\mypp
\bfc_T(s)+\bfq_T^{-}(\xi,s)\Bigr\}=o(1), \qquad \xi\to-\infty.
\end{equation}
Note that the expression in the curly brackets in
(\ref{eq:l'Hopital*}) has a finite limit as $\xi\to-\infty$, which
then must vanish in order to extinguish the multiplier
$\rme^{-\kappa\xi}\to\infty$, that is,
\begin{equation}\label{KK-*}
(I-\kappa^{-1}\myn A_T)\mypp\bfc_T(s)=-\bfq_T^{-}(-\infty,s).
\end{equation}
Conversely, condition (\ref{KK-*}) implies the limit
(\ref{eq:l'Hopital*}), since, by the l'H\^{o}pital rule, we have
\begin{equation*}
\frac{\bfq_T^{-}(\xi,s)-\bfq_T^{-}(-\infty,s)}{\rme^{\kappa\xi}}
\sim (I-\kappa^{-1}\myn A_T)\,\frac{f_0(cT\xi) \mypp\bbfw_T(\xi,s)}
{\kappa}=o(1), \qquad \xi\to-\infty.
\end{equation*}

Analogous considerations applied to (\ref{eqn:calB-plus*}) lead to
the integral equation
\begin{equation}\label{intsysu+}
\bfdelta_T(\xi,s)=\rme^{\xi\tilde{A}_T}\myp
\tildebfc_T+\int_0^\xi\rme^{(\xi-y)\myp \tilde{A}_T}
f_0(cTy)\mypp\bbfw_T(y,s)\,\rmd{y},\qquad \xi\ge0,
\end{equation}
with $\tildebfc_T\equiv\tildebfc_T(s)=\lim_{\xi\to 0+}
\bfdelta_T(\xi,s)$, which, similarly to (\ref{KK-*}), implies the
condition
\begin{equation}\label{KK+*}
(I+\tilde{\kappa}^{-1}\tilde{A}_T)\mypp\tilde{\bfc}_T(s)
=-\tilde{\bfq}_T^{+}(+\infty,s),
\end{equation}
where $\tilde\kappa=\kappa(\tilde s)$ and
\begin{equation}\label{eq:tilde-int_pm*}
\tilde{\bfq}_T^{\pm}(\xi,s):=
(I\pm\tilde{\kappa}^{-1}\tilde{A}_T)\int_0^\xi
\rme^{\mp\tilde{\kappa} y} f_0(cTy)\mypp\bbfw_T(y,s)\,\rmd{y},\qquad
\xi\ge0.
\end{equation}
Moreover, since the function $\bfdelta_T(\cdot,s)$ is continuous at
$\xi=0$, from formulas (\ref{intsysu-}) and (\ref{intsysu+}) we see
that $\bfc_T(s)=\tildebfc_T(s)$. Using this and subtracting
(\ref{KK-*}) from (\ref{KK+*}), we obtain
\begin{equation}\label{eq:c=}
\bfc_T(s)= \bigl(\kappa^{-1}\myn
A_T+\tilde{\kappa}^{\myp-1}\tilde{A}_T\bigr)^{-1}
\Bigl[\bfq_T^{-}(-\infty,s)-\tilde{\bfq}_T^{+}(+\infty,s)\Bigr].
\end{equation}

Evaluating the matrix inverse in (\ref{eq:c=}) is facilitated by
introducing the matrices (suggested by formulas (\ref{Ato1}))
\begin{equation}\label{eq:K-m}
K:=\left(
\begin{array}{rc}
1& 1\\[.3pc]
\!-1& 1
\end{array}
\right),\qquad K^{-1}=\frac{1}{2}\left(
\begin{array}{cr}
1\!&\! -1\\[.3pc]
1\!&\! 1
\end{array}
\right)
\end{equation}
and observing that
$$
K^{-1}\myn A_T\myp K=\kappa \left(\begin{array}{cc}
0 & s/\kappa\\[.3pc]
\kappa/s & 0
\end{array}\right).
$$
This gives
\begin{equation}\label{eq:R}
K^{-1}(\kappa^{-1}\myn A_T+\tilde{\kappa}^{-1}\tilde{A}_T)\myp K
=(s\tilde\kappa+\tilde{s}\kappa)\,\calR_T^{-1}(s), \qquad
\calR_T(s):=\left(
\begin{array}{cc}
0& s\tilde s\\[.3pc]
\kappa\tilde\kappa& 0
\end{array}
\right),
\end{equation}
and, returning to (\ref{eq:c=}), we finally get
\begin{equation}\label{eq:c=!}
\bfc_T(s)=(s\tilde\kappa+\tilde s\kappa)^{-1} K\mypp \calR_T(s)
K^{-1}\Bigl[\bfq_T^{-}(-\infty,s)-\tilde{\bfq}_T^{+}(+\infty,s)\Bigr].
\end{equation}

In view of Theorem \ref{th3} and according to (\ref{w-w*}), to
complete the proof of Theorem \ref{th4} we have to check that if
$\xi\sqrt{\lambda T}\to a\in\RR$ as $T\to\infty$ then
$\bfdelta_T(\xi,s)\to{\bf 0}$. To this end suppose, for instance,
that $\xi\le0$ and $a\le 0$ (the mirror case $\xi\ge0$, $a\ge 0$ is
considered similarly). Note that, as $T\to\infty$,
\begin{equation}\label{eq:xi*kappa*}
\kappa\sim \sqrt{2s\lambda T}, \qquad \kappa\mypp\xi=
\sqrt{2s}\,a+o(1),\qquad \kappa^{-1}\myn A_T= \lambda T
\kappa^{-1}J_1+O(\kappa^{-1}).
\end{equation}

Recall that the vectors $\bfq_T^{\pm}(\xi,s)$,
$\tilde{\bfq}_T^{+}(\xi,s)$ are defined in (\ref{eq:int_pm*}),
(\ref{eq:tilde-int_pm*}), respectively.

\begin{lemma}\label{lm:pmI*}
For each $s>0$\myp\textup{,} 
$\bfq_T^{-}(-\infty,s)=o(1)$ and
$\tilde{\bfq}_T^{+}(+\infty,s)=o(1)$ as $T\to\infty$.
\end{lemma}
\proof Both $\bfq_T^{-}$ and $\tilde{\bfq}_T^{+}$ are considered
similarly. For instance, using (\ref{eq:xi*kappa*}) and making the
change of variable $z=\kappa\myp y$, we have
\begin{equation}\label{eq:I->0}
|\bfq_T^{-}(-\infty,s)|= O(1) \int_{-\infty}^0 \rme^{z}
\bigl|f_0(cT\kappa^{-1}z)\bigr|\,\rmd{z}=o(1),\qquad T\to\infty,
\end{equation}
since, by the assumption of Theorem \ref{th4}, $cT\kappa^{-1}
\sim(2s)^{-1/2} \sqrt{c^2 T/\lambda}\to\infty$, hence
$f_0(cT\kappa^{-1}z)\to0$ for each $z<0$, and we can apply
Lebesgue's dominated convergence theorem.
\endproof

\begin{lemma}\label{lm:Ixi*}
As $T\to\infty$\myp\textup{,} if\/ $a_T:=\xi\sqrt{\lambda T}\to
a\in\RR$ then\myp\textup{,} for each $s>0$\myp\textup{,}
$\bfq_T^{\pm}(\xi,s)\to0$.
\end{lemma}
\proof By the substitution $y=\xi z$ and with the help of asymptotic
relations (\ref{eq:xi*kappa*}), we have
\begin{align}
\notag \bfq_T^{\pm}(\xi,s)&= \pm\bigl(\lambda
T\kappa^{-1}J_1+O(1)\bigr)\,\xi\int_0^1 \rme^{\mp\kappa\xi
z} f_0(cT\xi z)\mypp\bbfw_T(\xi z,s)\,\rmd{z}\\
\label{eq:a+o(1)*}
&=O(1)\,a_T\int_0^1 \bigl|f_0\bigl(z\myp a_T\sqrt{c^2T/\lambda}\,
\bigr)\bigr|\,\rmd{z}.
\end{align}
Now, if $a_T\to a=0$ then the right-hand side of (\ref{eq:a+o(1)*})
vanishes in the limit as $T\to\infty$, since the function $f_0$ is
bounded. If $a_T\to a\ne0$ then, similarly to the proof of Lemma
\ref{lm:pmI*}, the integral in (\ref{eq:a+o(1)*}) tends to zero
thanks to Lebesgue's dominated convergence theorem.
\endproof

Let us now return to equation (\ref{eq:D12*}). Using the identity
(\ref{KK-*}) and regrouping, we have
\begin{align}\notag
\bfdelta_T(\xi,s)&=\rme^{\kappa\xi}\mypp\kappa^{-1}\myn
A_T\mypp\bfc_T(s)-\cosh(\kappa\xi)\,\bfq_T^{-}(-\infty,s)+\frac{1}{2}\,
\rme^{\kappa\xi}\,\bfq_T^{+}(\xi,s)
+\frac{1}{2}\,\rme^{-\kappa\xi}\,\bfq_T^{-}(\xi,s)\\[.2pc]
\label{eq:delta} &= O(1)\mypp\kappa^{-1}\myn
A_T\mypp\bfc_T(s)+o(1),\qquad T\to\infty,
\end{align}
according to the second asymptotic relation in (\ref{eq:xi*kappa*})
and Lemmas \ref{lm:pmI*} and \ref{lm:Ixi*}. Further, substituting
the expression (\ref{eq:c=!}) for $\bfc_T$ and using Lemma
\ref{lm:pmI*} and the last relation in (\ref{eq:xi*kappa*}), we
obtain
\begin{equation}\label{eq:Ac}
\kappa^{-1}\myn
A_T\mypp\bfc_T=\frac{1}{\kappa\myp(s\tilde\kappa+\tilde
s\kappa)}\,\bigl(\lambda T J_1+O(1)\bigr) K\mypp \calR_T
K^{-1}o(1),\qquad T\to\infty.
\end{equation}
In turn, using the expressions (\ref{eq:K-m}) for the matrices $K$
and $K^{-1}$ and recalling the definition of the matrix $\calR_T$
given in (\ref{eq:R}), it is easy to calculate
\begin{equation}\label{eq:KRK}
K\mypp \calR_T K^{-1}=\frac{\kappa \tilde{\kappa}}{2}\,
J_1+O(1),\qquad T\to\infty.
\end{equation}
Finally, combining (\ref{eq:Ac}) and (\ref{eq:KRK}) and noting that
$J_1^2=0$ (see (\ref{eq:J12})), we have $\kappa^{-1}\myn
A_T\mypp\bfc_T=o(1)$ and hence, from (\ref{eq:delta}),
$\bfdelta_T(\xi,s)=o(1)$ as required. This completes the proof of
Theorem \ref{th4}.

\section{Concluding remarks}\label{sec:Concl}
We performed computer simulations to illustrate numerically the
convergence to the arcsine law, as stated by Theorems \ref{th3} and
\ref{th4}, for the occupation time functional
$\eta_T^{\pm}(0;f)=T^{-1}\!\int _0^T f(X_t^{\pm})\,\rmd{t}$ with
various probing functions $f$. The simulation algorithm is easily
implemented
by virtue of the obvious decomposition
$$
T\eta_T^{\pm}(0,f)\equiv\int_{0}^{T}f(X_{t})\,\rmd{t}=
\sum_{i=0}^{n-1}\int_{0}^{\tau_{i+1}}\!f\bigl(X_{\sigma_i}\mynn+(-1)^ict\bigr)\,\rmd{t}
+\int_{0}^{T-\sigma_{n}}\!f\bigl(X_{\sigma_n}\mynn+(-1)^nct\bigr)\,\rmd{t}\myp,
$$
where $(\tau_i)$ is a sequence of independent random times with
exponential distribution each (with parameter $\lambda$), and
$\sigma_i:=\tau_1+\cdots+\tau_i$ are the successive reversal times
of the telegraph motion; the threshold value $n$ is determined by
the condition $\sigma_n\le T<\sigma_{n}+\tau_{n+1}$.

\begin{figure}[h]
\centering
\includegraphics[width=.41\textwidth]{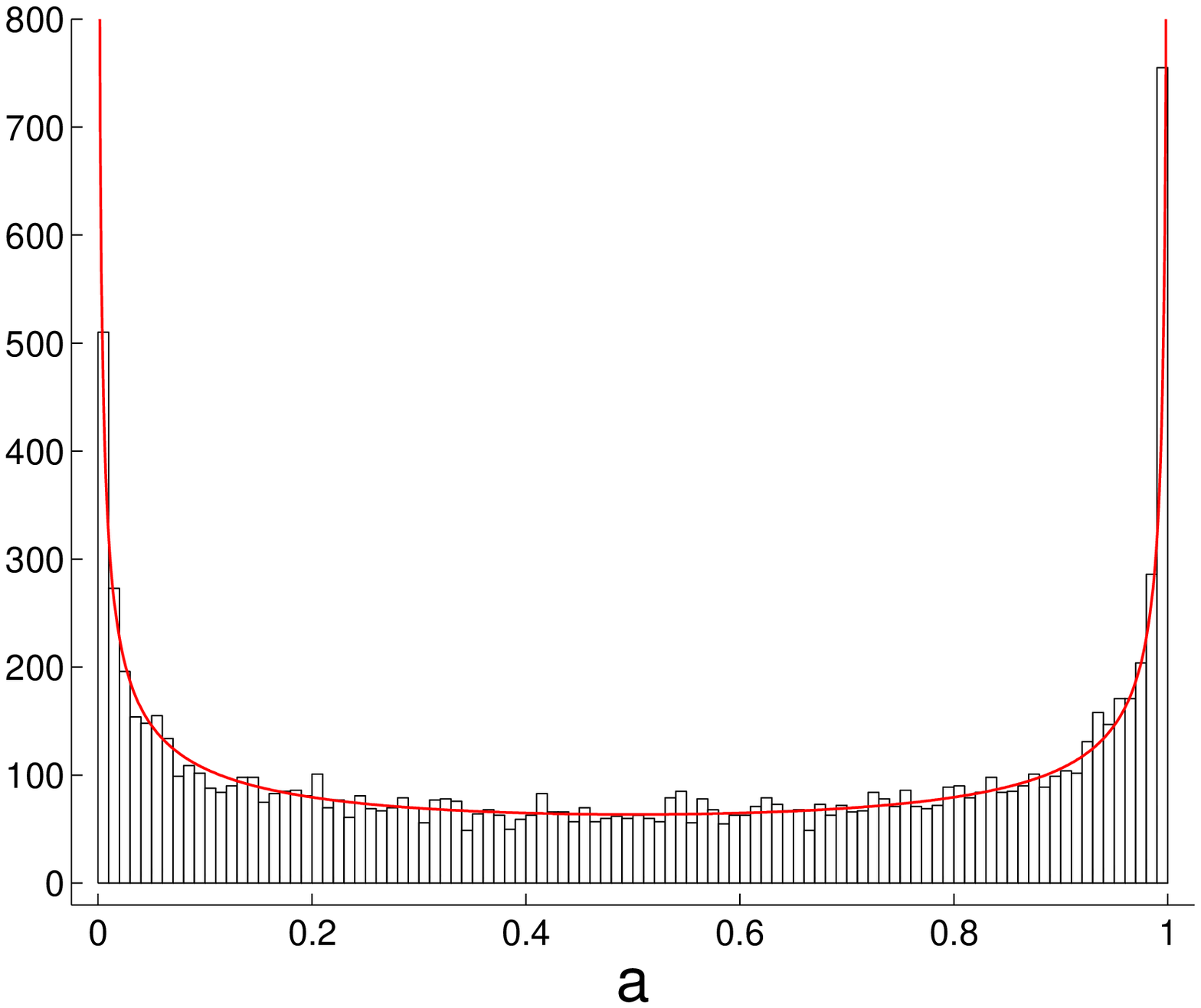}
\hspace{3pc}
\includegraphics[width=.41\textwidth]{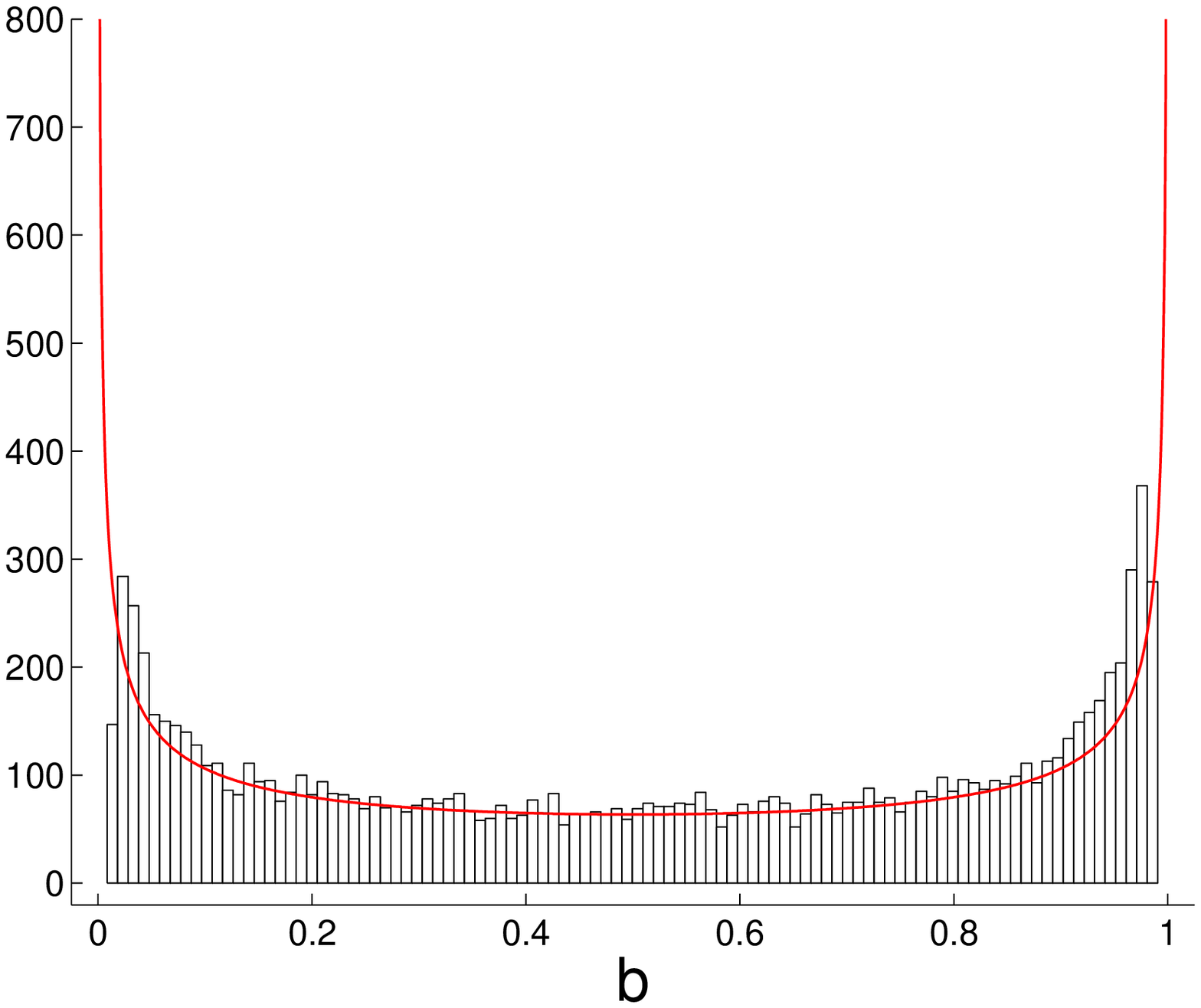}
\caption{\,Histograms for the occupation time functional
$\eta_T^{+}(0;f)$ with (a) the Heaviside step function $f=H$ and (b)
the function $f(x)=\pi^{\myp-1}\myn\arctan x+\frac{1}{2}$\myp. The
parameters of the telegraph process $X_t^{+}$ are standardized to
$c=1$ and $\lambda=1$. Both histograms are obtained with
$N=10,\myn000$ simulations, each over the observation time
$T=1\mypp000$. The length of each box on the histogram is
$\Delta=0.01$. The red solid curve represents the scaled arcsine
density (i.e., multiplied by $N\Delta=100$).}\label{fig1}
\end{figure}

\begin{figure}[h]
\centering
\includegraphics[width=.41\textwidth]{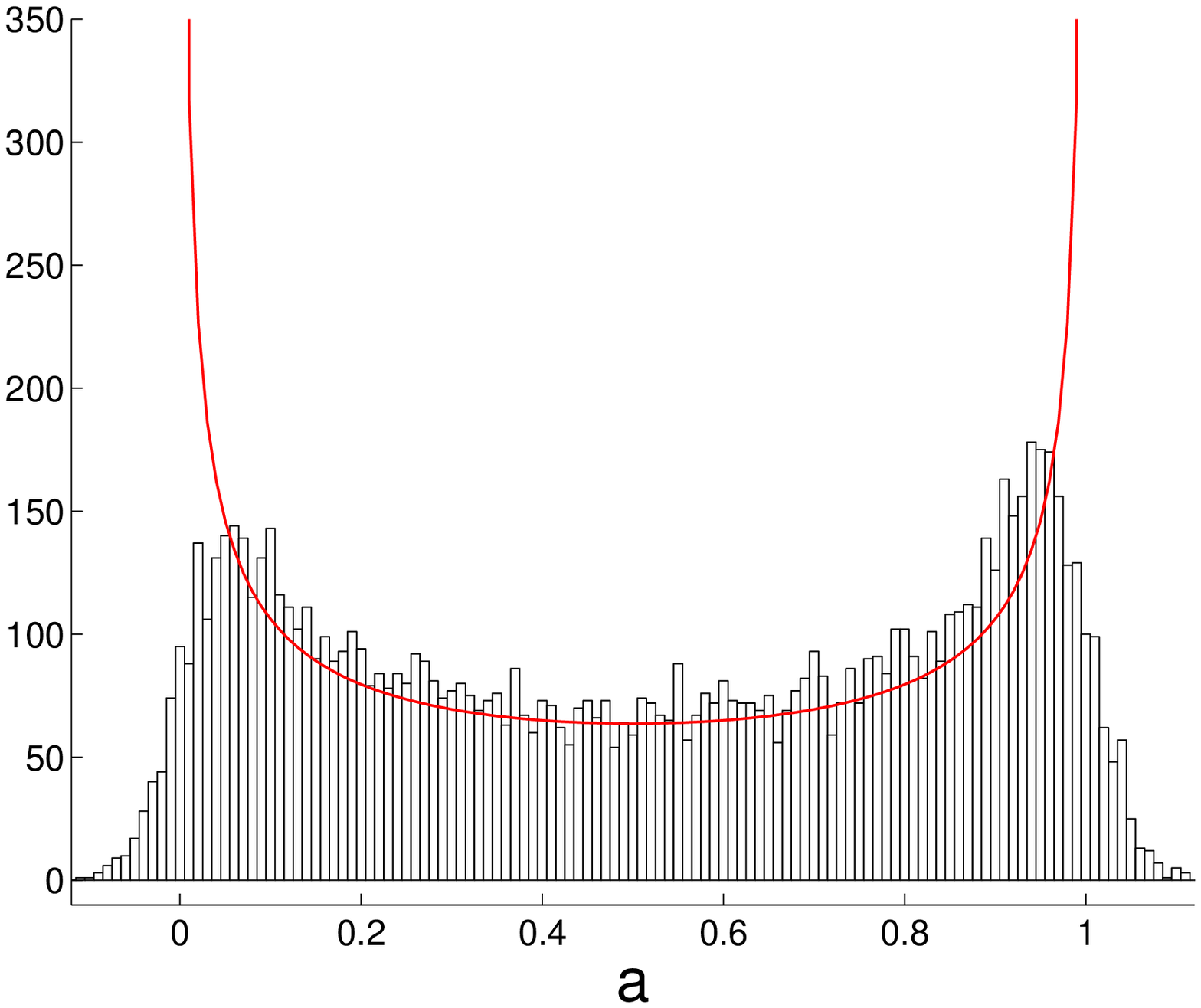}
\hspace{3pc}
\includegraphics[width=.41\textwidth]{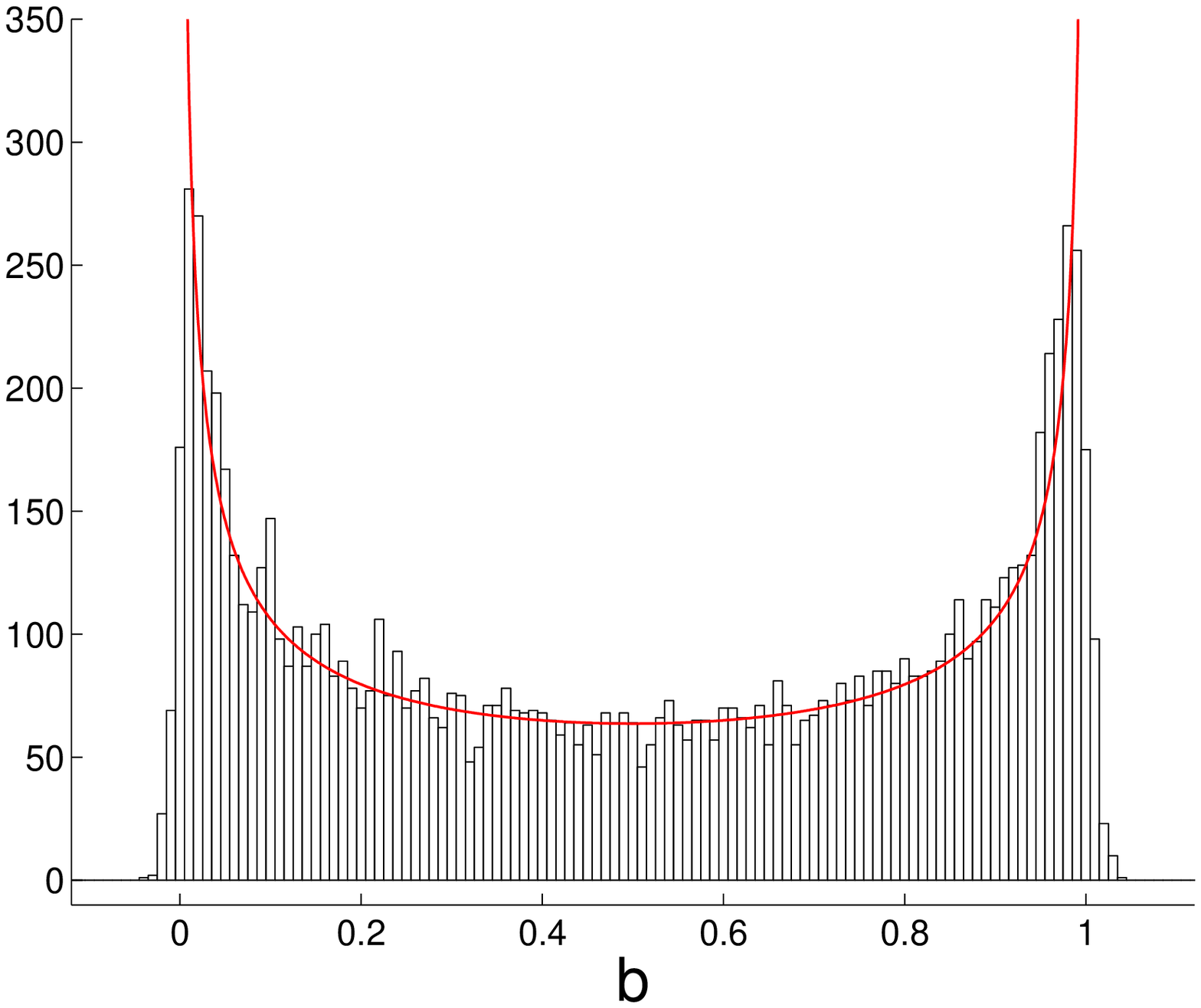}
\caption{\,Histograms for the functional $\eta_T^{+}(0;f)$ with the
probing function $f(x)=\pi^{\myp-1}\mynn\arctan x+\cos
x+\frac{1}{2}$\myp. The parameters of the telegraph process are as
in Figure \ref{fig1}, with the same number of runs $N=10,\myn000$
and the observation time (a) $T=1\mypp000$ or (b) $T=10\mypp000$.
Compare with Figure \ref{fig1} and note the improved quality of fit
to the hypothetical arcsine distribution (red curve) on the right
plot as compared to the left one.}\label{fig2}
\end{figure}

Throughout the simulations, we used the standardized parameters
$c=1$, \,$\lambda=1$, and plotted histograms of the sample values of
$\eta_T^{\pm}(0;f)$ based on $N=10,\myn000$ runs of the telegraph
process. To be specific, we simulated the plus-version of the
process, $X_t^{+}$ (i.e., with positive initial velocity), leading
to histograms slightly skewed to the right, especially at moderate
times $T$. No formal goodness-of-fit tests were applied, but the
histograms in Figures \ref{fig1} and \ref{fig2} below clearly
demonstrate the developing $U$-shape characteristic of the arcsine
distribution, however with the speed of such a convergence
apparently depending on the function $f$ involved (and, of course,
on the observation time $T$ used).

We start with the ``canonical'' case where the Heaviside function
$H(x)=\bfOne_{(0,\infty)}(x)$ plays the role of the probing function
$f$. Simulated values of $\eta_T^{+}(0;H)$ were obtained over the
observation time $T=1\mypp000$. The histogram plotted in Figure
\ref{fig1}{a} shows a very good fit of the data to the theoretical
arcsine density (rescaled according with the chosen representation
of the histogram). As already mentioned, the noticeable difference
between the highest columns at the left and right edges may be
attributed to asymmetry of the process $X_t^{+}$. More precisely,
the proportion of the sample values of $\eta_T^{+}(0;H)$ falling,
say, in the first box, $\Delta_1$ (from $0$ to $0.01$) and the last
box, $\Delta_{100}$ (from $0.99$ to $1$) is given by 510 and 750,
respectively, yielding the relative frequencies
$510/10,\myn000=0.051$ and $750/10,\myn000=0.075$. The corresponding
limiting probabilities, computed from the arcsine distribution
(\ref{eq:arcsine}), equal
$0.064$ for both $\Delta_1$ and $\Delta_{100}$ (here and below, we
give numerical values to two significant figures). This discrepancy
can be quantified using the exact theoretical distribution of
$\eta_T^{+}(0;H)$ obtained in Theorem \ref{th1} (see formula
(\ref{eta0a}) with $T=1\mypp000$), giving
the probability $0.052$ for $\Delta_1$ and $0.077$ for
$\Delta_{100}$, where the latter includes the atom
$2\myp\varphi_T(1)=0.025$. For comparison, with a tenfold
observation time
$T=10\mypp000$, these probabilities become 
$0.060$ and
$0.068$, respectively, with the atom much reduced,
$0.008$. It is also worth mentioning that, as indicated by these
results, the fit with the limiting arcsine distribution would be
much better for the ``symmetric'' version $\eta_T(0;H)$
corresponding to the telegraph process $X_t$ (see (\ref{VX})).

The long-term prediction contained in a more general Theorem
\ref{th4} was verified by computer simulations for the functional
$\eta_T^{+}(0;f)$ with the probing function
$f(x)=\pi^{-1}\mynn\arctan x+\frac12$\myp. The new histogram plot
(see Figure \ref{fig1}b), obtained with the same values of $c$,
$\lambda$, $T$ and $N$, is qualitatively similar to that on Figure
\ref{fig1}a, including a small right bias, but convergence to the
arcsine distribution becomes slower, apparently due to additional
time needed for the process to explore the limiting values $f_{\pm}$
of the function $f$ at $\pm\infty$, which eventually determine the
distributional limit.

Incidentally, this observation helps to understand the difference
between the sets of hypotheses in Theorems \ref{th3} and
\ref{th4}\myp; indeed, the additional condition of Theorem
\ref{th4}, requiring that $c^2T\myn/\lambda\to\infty$ as
$T\to\infty$, guarantees a sufficient mobility of the telegraph
process needed to gauge the limits $f_{\pm}$ available only at
remote distances from the origin. In contrast, if the function $f$
is reduced to the Heaviside step function $H$, the limiting values
$H_{-}\mynn=0$, $H_{+}\mynn=1$ are encountered by the process
straight away, so no extra mobility is needed.

Let us point out that the asymptotic conditions imposed in Theorem
\ref{th4} on the function $f$ involved in the occupation functional
$\eta_T^{\pm}(x;f)$  are rather strong, assuming the existence of
the limits $\lim_{x\to\pm\infty}f(x)=f_{\pm}$ \,(with $f_{-}\ne
f_{+}$)\,. This is in contrast with the paper by Khasminskii
\cite{Khas} mentioned in the Introduction, where the function $f$ is
only assumed to be Ces\`aro $(C, 1)$-summable, i.e., subject to a
weaker condition
$\lim_{x\to\pm\infty} \,x^{-1}\!\int_0^x f(u)\,\rmd{u}=f_{\pm}$
\,(cf.\ (\ref{Cesaro})). Unfortunately, we were unable to reach the
same level of generality. In particular, our proofs of formulas
(\ref{KK-*}), (\ref{KK+*}) and the key Lemmas \ref{lm:pmI*} and
\ref{lm:Ixi*} (see Section \ref{sec:th4}) are heavily based on the
existence of the limits $f_{\pm}$.

However, we conjecture that Theorem \ref{th4} does hold under the
weaker condition of Ces\`aro $(C, 1)$-summability of the probing
function $f$. To verify this claim numerically, we carried out
computer simulations for the distribution of $\eta_T^{+}(0;f)$ with
$f(x)=\pi^{\myp-1}\mynn \arctan x+\cos x+\frac{1}{2}$\mypp. Figure
\ref{fig2}\myp{a} shows the simulated histogram with the old values
$T=1\mypp000$ and $N=10,\myn000$, which reveals a bimodal
distribution but not quite well fit to the hypothetical arcsine
limit; in particular, there are noticeable ``parasite'' shoulders
outside the interval $[0,1]$, which are indeed possible because the
function $f$ may take values less than $0$ and bigger than $1$.
However, the fit with the arcsine shape significantly improves under
longer observations, $T=10\mypp000$ (see Figure \ref{fig2}\myp{b}).
In particular, the high modes at the edges are better pronounced,
while the shoulders outside $[0,1]$ are considerably reduced.

\appendix
\section{Probabilistic proof of Theorem \ref{th2}}\label{sec:A1}

Let us recall some information related to the first-passage problem
for the telegraph process $X^{\pm}_t$. For $x<0$, let
$\mathfrak{T}^{\pm}_{-x}:=\min\{t\ge0: X^{\pm}_t=-x\}$ (with the
convention that $\inf\myp\emptyset:=+\infty$) be the hitting time of
point $-x>0$ by the process $X^{\pm}_t$ (starting from the origin,
$X^{\pm}_0=0$). If we set $T_0:=(-x)/c$, then the distribution of
$\mathfrak{T}_{-x}^{\pm}$ is concentrated on $[T_0,\infty)$ and is
given by (see \cite[\S\mypp0.5, pp.\,12--13]{Pinsky} and also
\cite[Theorem 4.1, p.\,18]{Or95})
\begin{equation}\label{eq:T+}
\PP\{\mathfrak{T}^{+}_{-x}\in\rmd{t}\}= \rme^{-\lambda T_0}
\myp\delta_{T_0}(\rmd{t})+Q_{-x}^+(t)\,\rmd{t},\qquad
\PP\{\mathfrak{T}^{-}_{-x}\in\rmd{t}\}= Q_{-x}^-(t)\,\rmd{t},
\end{equation}
where the densities $Q_{-x}^{\pm}$ are defined exactly by equations
(\ref{eq:Q+}), (\ref{eq:Q-}).

Consider the two-dimensional Markov process $(X_t^{\pm},V_t^{\pm})$,
where $X_t^{\pm}$ is the (conditional) telegraph process
(\ref{VX-pm}) (i.e., with the initial velocity $V_0=\pm c$,
respectively), and $V_t^{\pm}=\rmd X^{\pm}_t/\rmd{t}=\pm c\myp
(-1)^{N_t}$ is the corresponding velocity process driven by an
dependent Poisson process $N_t$ which determines the reverse
instants of the motion $X^{\pm}_t$ (see (\ref{VX-pm})). It is
obvious that $\mathfrak{T}^{\pm}_{-x}$ is a stopping time for the
process $(X^{\pm}_t,V_t^{\pm})$. Also note that
$V^{\pm}_t\bigr|_{t=\mathfrak{T}^{\pm}_{-x}}\!=+c$ \,(a.s), since
the first passage through point $-x>0$ by the process $X^{\pm}_t$,
starting from the origin, with probability $1$ can only occur from
left to right, that is, with \textit{positive} velocity. Hence,
conditioning on the hitting time of the origin starting from $x<0$
(which, of course, has the same distribution as
$\mathfrak{T}^{\pm}_{-x}$) and using the strong Markov property of
the joint process $(X_t^{\pm},V_t^{\pm})$, we have, for each
$y\in[0,1-T_0/T]$,
\begin{align}
\notag \PP\{\myp{}&\eta_T^{\pm}(x)\in
\rmd{y}\}=\PP\{\mathfrak{T}^{\pm}_{-x}>T\}\myp\delta_{0}(\rmd{y}) +
\EE\bigl[\PP\{\eta^{\pm}_T(x)\in\rmd{y}, \,T_0\le
\mathfrak{T}_{-x}^{\pm}\le T\,|\,\mathfrak{T}_{-x}^{\pm}\}\bigr]\\
\label{eq:Y-new}
&=\left(\int_{T}^\infty\PP\{\mathfrak{T}^{\pm}_{-x}\in\rmd{u}\}\right)
\delta_{0}(\rmd{y}) + \int_{T_0}^{(1-y)\myp
T}\PP\{\mathfrak{T}^{\pm}_{-x}\in\rmd{u}\}
\,\PP\{(1-u/T)\,\eta^+_{T-u}(0)\in\rmd{y}\}.
\end{align}
Here, the first integral represents the case where the telegraph
process $X^{\pm}_t$ does not reach the origin before time $T$ and,
therefore, never enters the positive half-line (thus contributing to
the atom $\delta_{0}(\rmd{y})$), while the second integral (where
integration is taken with respect to $\rmd{u}$) accounts for the
first passage event (at time instant $u\in[T_0,(1-y)T]$), so that
the telegraph process, restarted from the origin (with the initial
velocity $+c$), has to spend on the positive half-line the required
time $T\myp\rmd{y}$ during the remaining travel time $T-u$.

In view of (\ref{eq:T+}) together with (\ref{eq:Q+}) and
(\ref{eq:Q-}), and due to equation (\ref{eta0a}) which provides the
distribution of $\eta_{T-u}^+(0)$, formula (\ref{eq:Y-new})
furnishes an explicit representation of the distribution of
$\eta_T^{\pm}(x)$. More explicitly, on account of the atom in
(\ref{eq:T+}), the right-hand side of (\ref{eq:Y-new}) specializes
to
\begin{equation}\label{eq:Y-*}
\left(\int_T^\infty
Q^{\pm}_{-x}(u)\,\rmd{u}\right)\delta_{0}(\rmd{y})+
\mu^{\pm}_T(\rmd{y}) +\int_{T_0}^{(1-y)\myp T}
\!Q^{\pm}_{-x}(u)\;\PP\left\{\eta^+_{T-u}(0)\in
\frac{\rmd{y}}{1-u/T}\right\}\rmd{u}\myp,
\end{equation}
where $\mu^{-}_T(\rmd{y}):=0$ and
\begin{equation}\label{eq:mu+*}
\mu_T^{+}(\rmd{y}):= \rme^{-\lambda T_0}
\,\PP\left\{\eta^+_{T-T_0}(0)\in \frac{\rmd{y}}{1-T_0/T}\right\}.
\end{equation}
Using (\ref{eta0a}), for any $u\in[T_0,(1-y)\myp T]$ we have
\begin{equation}\label{eq:atom}
\PP\left\{\eta^+_{T-u}(0)\in
\frac{\rmd{y}}{1-u/T}\right\}=2\varphi_{T-u}(1)\,\delta_{1-u/T}(\rmd{y})
+\psi_{T-u}\left(\frac{y}{1-u/T}\right)\frac{\rmd{y}}{1-u/T}\,.
\end{equation}
Substituting (\ref{eq:atom}) (with $u=T_0$) into (\ref{eq:mu+*})
readily gives (\ref{eq:mu}), while the last term on the right-hand
side of (\ref{eq:Y-*}) is reduced to (cf.\ (\ref{eq:Y-}))
\begin{equation*}
2\myp T\myp Q^{\pm}_{-x}((1-y)T)\,\varphi_{T}(y)\,\rmd{y}+
\int_{T_0}^{(1-y)\myp T} Q^{\pm}_{-x}(u)\left(\psi_{T-u}
\left(\frac{y}{1-u/T}\right)\frac{\rmd{y}}{1-u/T}\right)\rmd{u}\myp,
\end{equation*}
where the contribution of the atom $\delta_{1-u/T}(\rmd{y})$ from
(\ref{eq:atom}) is easily computed via the obvious symbolic formula
$\delta_{1-u/T}(\rmd{y})\,\rmd{u}=T\delta_{(1-y)T}(\rmd{u})\,\rmd{y}$\myp.
Indeed, for any test functions $F(y)$ and $G(u)$ we have, by
changing the order of integration,
\begin{align*}
\int_0^{1-T_0/T}\!F(y)&\int_{T_0}^{(1-y)\myp T}
\!G(u)\,\delta_{1-u/T}(\rmd{y})\,\rmd{u}=\int_{T_0}^{T}
G(u)\,\rmd{u}\int_0^{1-u/T}F(y)\,\delta_{1-u/T}(\rmd{y})\\
&=\int_{T_0}^{T} G(u)\,F(1-u/T)\,\rmd{u} =T\int_{0}^{1-T_0/T} F(y)
\left(\int_{T_0}^{(1-y)T} G(u)\,\delta_{(1-y)\myp
T}(\rmd{u})\right)\rmd{y}\myp.
\end{align*}

\section{Probabilistic proof of Theorem \textup{\ref{th3}}}\label{sec:A2}
Making the substitution $t=Tu$ and using that $H(\alpha x)\equiv
H(x)$ for any $\alpha>0$, we can rewrite formula (\ref{OT0}) as
\begin{equation}\label{eq:rescaled}
\eta_T^\pm(x) = \int_0^1 H(x+X^\pm_{Tu})\,\rmd{u} = \int_0^1
H(\tilde{x}+\tilde{X}^\pm_{u})\,\rmd{u},
\end{equation}
where $\tilde{x}:=(c^2T/\lambda)^{-1/2}x$, and
$\tilde{X}_u^\pm:=(c^2T/\lambda)^{-1/2}X_{Tu}^\pm$ ($u\ge 0$) is
another telegraph process with rescaled parameters
$\tilde{\lambda}:=\lambda T\to\infty$, \,$\tilde{c}:=(\lambda
T)^{1/2}\to\infty$ \,($T\to\infty$). By Theorem \ref{Kac}, the
process $(\tilde X_u^\pm,\,0\le u\le1)$ converges weakly to a
standard Brownian motion $(B_u,\,0\le u\le1)$. Hence, if
$\tilde{x}\to a$ as $T\to\infty$ (cf.\ the hypotheses of Theorem
\ref{th3}) then from (\ref{eq:rescaled}) we immediately obtain the
convergence in distribution, as $T\to\infty$,
\begin{equation*}
\eta_T^\pm(x)\stackrel{\!d\,}{\to} \fL_1(a):=\int_0^1
H(a+B_{u})\,\rmd{u}.
\end{equation*}

According to (\ref{eq:h_T}) and (\ref{Levy}), the random variable
$\fL_1(0)$ has the arcsine distribution, which proves Theorem
\ref{th3} for $a=0$. For $a<0$ (so that $-a>0$), let
$\tau_{-a}:=\min\{t\ge0: B_t=-a\}$ be the hitting time of the point
$-a$ by the Brownian motion $B_t$ starting from the origin
($B_0=0$). As is well known since P.~L\'evy's paper
\cite[Th\'eor\`eme~2, p.~294]{Levy} (see also \cite[\S\myp1.7,
p.~26]{Ito} or \cite[\S\myp{}VI.2\myp(e), pp.\,174--175]{Feller}),
the random variable $\tau_{-a}$ has probability density
$q_{-a}(\cdot)$ defined in (\ref{eq:q}). Note that $\tau_{-a}$ is a
stopping time (with respect to the natural filtration
$\calF_t:=\sigma\{B_s,\ 0\le s\le t\}$). Conditioning on $\tau_{-a}$
(when $a+B_{\tau_{-a}}=0$) and using the strong Markov property, we
obtain, for any $y\in[0,1]$,
\begin{align*}
\PP\{\myp{}\fL_1(a)\in
\rmd{y}\}&=\PP\{\tau_{-a}>1\}\,\delta_{0}(\rmd{y}) + \int_0^{1-y}
q_{-a}(u)\;\PP\{(1-u)\mypp\fL_{1-u}(0)\in \rmd{y}\}\,\rmd{u}\\
&=\left(\int_{1}^\infty
q_{-a}(u)\,\rmd{u}\right)\delta_{0}(\rmd{y})+ \left(\int_0^{1-y}
\frac{q_{-a}(u)}{1-u}\,p_{\mathrm{as}}\!\left(\frac{y}{1-u}\right)\rmd{u}\right)\rmd{y},
\end{align*}
which coincides with (\ref{eq:Y0}) (for $Y_{-a}$) in view of
(\ref{eq:m_a}) and (\ref{eq:f_a}). Finally, the case $a>0$ easily
follows by noting the obvious symmetry relation
$\fL_1(a)\stackrel{d}{=}1-\fL_1(-a)$ (cf.\ (\ref{eq:=d})).
\endproof

\begin{acknowledgements}
We gratefully acknowledge partial support by the London Mathematical
Society (through an LMS Scheme~2 Grant) during N.~Ratanov's visit to
the University of Leeds in June 2007, when part of this research was
done.
\end{acknowledgements}


\end{document}